\documentclass[fleqn]{article}
\usepackage{amssymb,latexsym,theorem, epic}

\newcommand{\C}{\mathcal{C}}

\newcommand{\PG}{\mathrm{PG}}

\newcommand{\F}{\mathbb{F}}

\newcommand{\E}{\mathcal{E}}

\newcommand{\T}{\mathcal{T}}

\newcommand{\X}{\mathcal{X}}

\newcommand{\ba}{\bar{a}}
\newcommand{\bb}{\bar{b}}
\newcommand{\bc}{\bar{c}}
\newcommand{\bd}{\bar{d}}
\newcommand{\be}{\bar{e}}

\newcommand{\bp}{\bar{p}}

\newtheorem{lemma}{Lemma}[section]
\newtheorem{prop}[lemma]{Proposition}
\newtheorem{theo}[lemma]{Theorem}
\newtheorem{theorem}{Theorem}
\newtheorem{co}[lemma]{Corollary}

\def\pr{\noindent{\bf Proof. }}
\def\eop{\hspace*{\fill}$\Box$}

\title{On elliptic ovoids and their rosettes in a classical generalized quadrangle of even order.}

\author{Ilaria Cardinali and N. S. Narasimha Sastry }
\date{}

\begin{document}
\maketitle
\begin{abstract}
Let $Q_0$ be the classical generalized quadrangle of order $q=2^n$ arising from a non-degenerate quadratic form in a $5$-dimensional vector space defined over a finite field of order $q$.
 We consider the rank two geometry  $\X$  having as points all the
elliptic ovoids of $Q_0$ and as lines the maximal pencils of elliptic ovoids
of $Q_0$ pairwise tangent at the same point. We first prove that $\X$ is isomorphic to a
2-fold quotient of the affine generalized quadrangle $Q\setminus Q_0$ where $Q$
is the classical $(q,q^2)$-generalized quadrangle admitting $Q_0$ as
a hyperplane. Then, we investigate the collinearity graph $\Gamma$ of $\X.$ In particular, we obtain a classification of the cliques of $\Gamma$ proving that they arise either from lines of $Q$ or subgeometries of $Q$ defined over $\F_2.$
\end{abstract}

\section{Introduction}\label{Introduction}
 We refer to the monograph~\cite{Payne-Thas} for the terminology and basics on finite generalized quadrangles.



We only recall that a finite generalized quadrangle $\mathcal{Q}$ is {\it classical} if its points and lines are points and lines of a projective space $\PG(d,q),$ the points of $\mathcal{Q}$ span $\PG(d,q)$ and all the lines of $\mathcal{Q}$ are full, namely if a projective line $l$ is a line of $\mathcal{Q}$ then all points of $l$ belong to $\mathcal{Q}.$  If this is the case then the points and the lines of $\mathcal{Q}$ are the points and lines of $\PG(d,q)$ that are singular and respectively totally singular for a suitable non-degenerate alternating, hermitian or quadratic form of the underlying vector space $V(d+1,q)$ of $\PG(d,q)$ (\cite[Chapter 4]{Payne-Thas}).

The generalized quadrangle of order $q$ arising from a non-degenerate symplectic form on $V(4,q)$ is denoted by $W(q).$  The generalized quadrangle of order $q$ arising from a non-degenerate quadratic form in $V(5,q)$ is denoted by $Q(4,q).$  The $(q,q^2)$-generalized quadrangle arising from a non-degenerate quadratic form of Witt index $2$ in $V(6,q)$ is denoted by $Q^-(5,q).$ The $(q,1)$-generalized quadrangle arising from a non-degenerate quadratic form of Witt index $2$ in $V(4,q)$ is denoted by $Q^+(3,q).$ It is also called a {\it grid}.
For any $q,$ the dual of $Q(4,q)$ is isomorphic to $W(q).$  Moreover, $Q(4,q)\cong W(q)$ if and only if $q$ is even.

An {\it ovoid} of an $(s,t)$-generalized quadrangle $\mathcal{Q}$ is a set $X$  of points of $\mathcal{Q}$ such that every line of $\mathcal{Q}$ meets $X$ in exactly one point (see~\cite{Payne-Thas}). Clearly, $|X|=st+1$ and no two points of $X$ are collinear. For a survey of ovoids of generalized quadrangles we refer the reader to~\cite{Payne-Thas-1994}.

In this paper we assume $q=2^n$  and we denote by $Q_0$ a classical generalized quadrangle
$Q(4,q)$ of order $q$ embedded in a $4$-dimensional projective space
$\PG(V_0),\,V_0=V(5,q).$ Let $n_0$ be the {\it nucleus} of $Q_0$, i.e.
$n_0$ is the unique non-singular point of $\PG(V_0)$ with the property
that the lines of $\PG(V_0)$ through $n_0$ are the tangents to $Q_0$ in $\PG(V_0)$.

Classification of ovoids in $Q_0$ is a fundamental open problem.
There
are two known classes of ovoids of $Q_0$, namely the class $\E$ of {\it
elliptic  ovoids} and the class $\T$ of {\it Tits-Suzuki ovoids.}
An ovoid $X$ of $Q_0$ is elliptic if and only if it spans a hyperplane $\langle X\rangle$ of $\PG(V_0)$ (whence $X=\langle X\rangle \cap Q_0$). Tits-Suzuki ovoids exist if and only if
$n$ is odd.

In this paper we will be concerned with a semipartial geometry arising from the class $\mathcal{E}$ of elliptic ovoids of $Q_0.$
Remind that any two distinct elliptic ovoids of $Q_0$ intersect in either a singleton or a non-degenerate conic.
If two elliptic ovoids meet at a single point $p$ we will say that they are {\it tangent} at $p$.



Following~\cite{BTVM-2002}, we define a {\it rosette $\rho_p$} of elliptic ovoids {\it based} at a point $p$ as a set of $q$ elliptic ovoids of $Q_0$ mutually intersecting at $p$ such that $\{X \setminus\{p\}\colon X \in \rho_p\}$ is a partition of the set of the points of $Q_0$ not collinear with $p.$ The point $p$ is called the {\it base point of $\rho_p$.} Clearly, a set of elliptic ovoids mutually tangent at a given point is a rosette if and only if it consists of exactly $q$ ovoids. It is also clear that a rosette with base point $p$ is a maximal pencil of elliptic ovoids of $Q_0$ mutually tangent at $p.$

A proof of the following proposition has been given in~\cite[Theorem 3.2]{Bagchi-Sastry-88} (see also~\cite[Example 1.4(d)]{Debroe-Thas-1978}). 

\begin{prop} \label{prop1-2}
Let $X_1$ and $X_2$ be two elliptic ovoids of $Q_0$ tangent at a point $p$. Then there exists a unique rosette $\rho_p$ with base point $p$ containing $X_1$ and $X_2.$
\end{prop}

\begin{prop} \label{prop1}
Let $\rho_p$ be a rosette of elliptic  ovoids of $Q_0$ based at the point $p.$  Then there exists a unique projective plane $\pi_p$ of $\PG(V_0)$ containing $p$ and tangent to each  $X\in \rho_p.$
\end{prop}
\pr
Let $X_i, X_j, X_k$ be three distinct elliptic ovoids contained in the rosette $\rho_p.$ Let $\pi_j:=\langle X_i \rangle \cap \langle X_j \rangle$ and $\pi_k:=\langle X_i \rangle \cap \langle X_k \rangle$ be the intersection planes of the $3$-dimensional subspace $\langle X_i \rangle$ spanned by $X_i$ with the  $3$-dimensional subspaces spanned by $X_j$ and $X_k$ respectively.
Since $\langle X \rangle \cap Q_0 =X$ for every $X\in \mathcal{E}$, we have $\pi_j\cap X_i=p= \pi_k\cap X_i.$

If $\pi_j\not= \pi_k$ then $\pi_j$ and $\pi_k$ are two distinct planes of $\langle X_i \rangle$ both tangent to $X_i$ at $p$. This is clearly impossible. Hence, $\pi_j=\pi_k$ and the proposition follows.\eop

\medskip

Let  $\X:=(\E,\mathcal{L})$ be the point-line geometry having the set $\E$ of elliptic ovoids of $Q_0$ as the
pointset and the set of rosettes of elliptic  ovoids of $Q_0$ as the lineset. An
ovoid $X$ is incident with a rosette $\rho_p$ if $X\in
\rho_p.$
It is known that $\X$ is a semi-partial geometry with parameters $(s,t,\alpha,\mu)
 = (q-1,q^2,2,2q(q-1))$ (see~\cite{Debroe-Thas-1978} and~\cite{DCVM}). A very straightforward proof of this claim  will also be given in Section~\ref{sec 3}.

\subsection{Main results}\label{main results-sec 1}
Let $Q\cong Q^-(5,q)$ be a classical $(q,q^2)$-generalized
quadrangle embedded in the $5$-dimensional projective space $\PG(V),$ where $V=V(6,q).$ We may assume that  $Q_0=Q\cap H_0$ for a suitable hyperplane  $H_0$ of $\PG(V).$  Thus, $Q_0$ is embedded in $H_0\cong\PG(4,q).$ Let $\perp$ denote the orthogonality relation defined by $Q$ in $\PG(V).$ The point $n_0:=H_0^\perp$ is the nucleus of the quadric $Q_0.$

Let $\widehat{\X}:=(\widehat{P},\widehat{L})$ be the point-line geometry
having the set $\widehat{P}$ of the points of $Q\setminus Q_0$ as the point-set and the set
$\widehat{L}:=\{l\setminus Q_0\colon l\,\,{\rm line\,\,of \,\,} Q {\,\,\rm not\,\,contained \,\,in\,\,Q_0} \}$ as line-set.

If $l$ is a line of $Q$ not contained in $Q_0$ then $l\cap Q_0$ is a point, henceforth denoted by $l^\infty.$ We denote by $\hat{l}$ the line of $\widehat{\X}$ corresponding to $l$ and we call $l^\infty$ the {\it point at infinity} of $\hat{l}.$ Thus,
\begin{equation}\label{l infinity}
\hat{l}=l\setminus \{l^\infty\},\,\ {\rm and}\,\, l=\hat{l}\cup \{l^\infty\}.
\end{equation}

The incidence relation in $\widehat{\X}$ is inherited from the incidence relation of $Q.$
It is easy to see that the incidence graph of $\widehat{\X}$ is connected. It has girth $8$, whence $\widehat{\X}$ has gonality $4$. The collinearity graph of $\widehat{\X}$ has diameter $3$ (see e.g.~\cite[8.4.1]{Pas-Diagram geometry}). \\

We recall that a {\it morphism} of rank 2 geometries is a morphism of their incidence graphs.
Given two rank 2 geometries $\mathcal{G}_1$ and $\mathcal{G}_2$, a
morphism $\varphi\colon \mathcal{G}_1\rightarrow \mathcal{G}_2$ is a {\it
covering} if, for any point $p$ of $\mathcal{G}_1$ and any line $l$ of
$\mathcal{G}_1$, $\varphi$ induces a  bijection from the pointset of  $l$ to the pointset of $\varphi (l)$ and from the
set of lines of $\mathcal{G}_1$ through $p$ to the set of lines of $\mathcal{G}_2$ through $\varphi(p).$  If the fibers of a cover have all the same
size $t$ then we will say that $\varphi$  is a {\it $t$-fold covering},
$\mathcal{G}_1 $ is called a {\it $t$-fold cover of $\mathcal{G}_2$} and $\mathcal{G}_2 $ is called a {\it $t$-fold quotient of $\mathcal{G}_1$}. Note that all coverings of connected geometries are surjective (see~\cite[8.3]{Pas-Diagram geometry}). In the terminology of
graph theory, a covering of geometries is a covering of the
incidence graphs of the given geometries. It also induces a covering of the collinearity graphs, satisfying the additional property that, if $a$, $b$, $c$ are different points of $\mathcal{G}_1$ with $b$ and $c$ collinear with $a$ and $\varphi(a)$, $\varphi(b)$, $\varphi(c)$ belonging to a common line of $\mathcal{G}_2$, then $a$, $b$, $c$ belong to the same line of $\mathcal{G}_1.$ We warn that, however, if $b$ and $c$ are collinear with $a$ and $\varphi(b)$ and  $\varphi(c)$ are mutually collinear but not on the same line as $\varphi(a)$, then $b$ and $c$ need not be collinear in $\mathcal{G}_1.$

In Section~\ref{sec 3} we shall prove the following:

\begin{theorem}\label{Main thm 1}
The geometry $\widehat{\X} $ is a $2$-fold cover of  the geometry $\X.$ Explicitly, let $\varphi_0\colon \widehat{\X}\rightarrow \X$ be the mapping defined as follows:
\[\varphi_0(x)=x^\perp\cap Q_0\,\,for\,\,every\,\,point\,\,x\in \widehat{P},\]
\[\varphi_0(\hat{l})=\rho_{p} \,\,where\,\, p=l^\infty,\,\,for\,\,every\,\,line\,\,\hat{l}\in \widehat{L}.\]
Then $\varphi_0$ is a $2$-fold covering.
\end{theorem}


The covering $\varphi_0$ defined as above will be called the {\it canonical} covering from $\widehat{\X}$ to $\X.$ As $Q$ can be recovered from $\widehat{\X}$ (see~\cite{Pralle}), every automorphism $\hat{\alpha}$ of $\widehat{\X}$ can be extended to a unique automorphism  $\alpha$ of $Q$ stabilizing $Q_0.$
Let $\alpha_0$ be the automorphism of $Q_0$ induced by $\alpha.$ Clearly, $\varphi_0 \circ \hat{\alpha}=\alpha_0\circ \varphi_0.$ It follows that the composition $\varphi:=\varphi_0\circ \hat{\alpha}$ is still a $2$-fold covering of $\widehat{\X}.$ Perhaps, all $2$-fold coverings from  $\widehat{\X}$ to $\X$ arise in this way, but we are not going to investigate this conjecture in this paper.

\medskip

In Section~\ref{collinearity graph} we determine the cliques of the collinearity graph
$\Gamma$ of $\X$.
A clique  of $\Gamma$ is called {\it linear} if it is a subset of a line of $\X$; it is said to be {\it non-linear} otherwise.
A clique of size $i$ is called an {\it $i$-clique}.

As we will prove in Section~\ref{collinearity graph} (Proposition~\ref{cliques-1}) a linear clique and a non-linear clique can have at most two vertices in common.\\

\begin{theorem}\label{Main thm 2}
All the following hold
\begin{itemize} \item[1)] Let $q=2^n$ with $n$ even. Then the maximal cliques of $\Gamma$ have size $q$ or $4.$ When $q>4$ cliques of size $q$ are linear while the maximal cliques of size $4$ are non-linear.
 \item[2)] Let $q=2^n >2$ with  $n$ odd. Then the maximal cliques of $\Gamma$ have size $q$ or $6.$ The cliques of size $q$ are linear cliques. The maximal cliques of size $6$ are non-linear.
\item[(3)] In any case, every non-linear $3$-clique can be extended to a non-linear
$4$-clique in $q+1$ ways.
\item[(4)] Let $q=2^n$ with $n$ odd. Then
\begin{itemize}
\item[(4.i)] Any non-linear $4$-clique can be extended to a
$5$-clique in two ways.
\item[(4.ii)] Any non-linear $5$-clique can be extended to a $6$-clique in only one way.
\end{itemize}
\end{itemize}

\end{theorem}

Note that the restrictions $q>4$ and $q>2$ in claim 1) and claim 2) of Theorem~\ref{Main thm 2} are essential. Indeed, when $q=2$ then $\Gamma$ is a complete graph with $6$ vertices. 
When $q=4,$ all maximal cliques of $\Gamma$ have size $4.$

\medskip
Denote by $\widehat{\Gamma}$ the collinearity graph of $\widehat{\X}.$

A {\it hexagon} $H=(a^1,b^1,c^1,a^2,b^2,c^2, a^1)$ of $\widehat{\X}$ is a $6$-circuit of
$\widehat{\Gamma}.$  Note that $H$ admits a bipartition by 3-sets $\{a^1,b^1,c^1\}$ and $\{a^2,b^2,c^2\}$ as well as a $3$-partition by $2$-sets $\{a^1, a^2\}$, $\{b^1,b^2\}$ and  $\{c^1,c^2\}$ such that vertices in the same pair are on a secant line to $Q$ through $n_0$. Vertices in the same pair are called {\it opposite}. We say that a hexagon is {\it centric} if all projective lines through opposite vertices are concurrent at the same point called the {\it center of the hexagon}.

\begin{picture}(310,112)(0,0)

\put(98,18){$\bullet$}
\put(96,8){$a^1$}

\put(100,21){\line(3,1){39}}

\dottedline{3}(100,19)(100,86)
\put(97,52){$\odot$}

\put(140,33){\line(0,1){41}}

\put(138,31){$\bullet$}
\put(145,31){$b^1$}

\dottedline{3}(140,34)(50,80)

\put(138,72){$\bullet$}
\put(145,72){$c^1$}

\put(100,88){\line(3,-1){39}}

\put(98,85){$\bullet$}
\put(96,93){$a^2$}

\put(60,34){\line(3,-1){39}}

\put(58,32){$\bullet$}
\put(48,32){$c^2$}
\dottedline{3}(60,36)(140,75)

\put(60,33){\line(0,1){41}}

\put(58,72){$\bullet$}
\put(48,72){$b^2$}

\put(60,74){\line(3,1){39}}
\put(86,-10){Figure 1}
\end{picture}
\bigskip
\medskip

A {\it cube} $C$ of $\widehat{\X}$ is an induced subgraph of
$\widehat{\Gamma}$ with $8$ vertices as in the following
picture, where adjacencies are represented by thick lines:


\begin{picture}(310,125)(0,0)
\put(0,18){$\bullet$}
\put(0,8){$a^1$}

\dottedline{3}(0,20)(100,100)

\dottedline{3}(0,80)(100,41)

\dottedline{3}(60,20)(40,100)

\dottedline{3}(60,80)(40,41)

\put(47,57){$\odot$}

\put(3,20){\line(1,0){57}}
\put(60,18){$\bullet$}
\put(60,8){$b^1$}
\put(1,20){\line(0,1){57}}
\put(61,20){\line(0,1){57}}
\put(0,88){$d^1$}
\put(60,85){$c^1$}
\put(3,80){\line(1,0){57}}
\put(0,77){$\bullet$}
\put(60,77){$\bullet$}
\put(2,20){\line(2,1){40}}
\put(62,20){\line(2,1){40}}
\put(2,79){\line(2,1){40}}
\put(62,79){\line(2,1){40}}
\put(40,38){$\bullet$}
\put(35,28){$c^2$}
\put(43,40){\line(1,0){57}}
\put(100,38){$\bullet$}
\put(100,28){$d^2$}
\put(41,40){\line(0,1){57}}
\put(101,40){\line(0,1){57}}
\put(40,108){$b^2$}
\put(103,108){$a^2$}
\put(43,100){\line(1,0){57}}
\put(40,97){$\bullet$}
\put(100,97){$\bullet$}
\put(49,-8){Figure 2}
 \end{picture}
\bigskip
\medskip

The graph $C$ admits a bipartition in $4$-sets $\{a^1,b^1,c^1,d^1\}$, $\{a^2,b^2,c^2,d^2\}$ as well as a $4$-partition in pairs $\{a^1, a^2\}$, $\{b^1,b^2\}$, $\{c^1,c^2\}$, $\{d^1,d^2\},$ two vertices in the same pair being called {\it opposite}. We say that a cube  is {\it centric} if all projective lines through opposite vertices are concurrent at the same point called the {\it center of the cube}.

Note that a hexagon (respectively, a cube) is the complement of a {\it ($3\times 2$)-grid graph} (respectively, a ($4\times 2$)-grid graph). The bipartition in two classes of size $3$ (respectively $4$) is one of the two families of lines of the grid and the $3$-partition ($4$-partition) in pairs is the other family. We exploit this observation to define {\it dodecades}. Consider a ($6\times 2$)-grid as follows:

\begin{picture}(310,100)(0,0)
\put(0,21){$\bullet$}
\put(0,8){$a^1$}
\put(3,24){\line(1,0){250}}

\put(53,21){$\bullet$}
\put(53,8){$b^1$}

\put(103,21){$\bullet$}
\put(103,8){$c^1$}

\put(153,21){$\bullet$}
\put(153,8){$d^1$}

\put(203,21){$\bullet$}
\put(203,8){$e^1$}

\put(253,21){$\bullet$}
\put(253,8){$f^1$}

\put(1,23){\line(0,1){50}}
\put(54,23){\line(0,1){50}}
\put(104,23){\line(0,1){50}}
\put(154,23){\line(0,1){50}}
\put(204,23){\line(0,1){50}}
\put(254,23){\line(0,1){50}}
\put(0,71){$\bullet$}
\put(0,78){$a^2$}
\put(53,71){$\bullet$}
\put(53,78){$b^2$}

\put(103,71){$\bullet$}
\put(103,78){$c^2$}

\put(153,71){$\bullet$}
\put(153,78){$d^2$}

\put(203,71){$\bullet$}
\put(203,78){$e^2$}

\put(253,71){$\bullet$}
\put(253,78){$f^2$}

\put(3,73){\line(1,0){250}}
\put(106,-8){Figure 3}
\end{picture}
\bigskip

A {\it dodecade} of $\widehat{\X}$ is an induced subgraph of
$\widehat{\Gamma}$ with $12$ vertices
\[D=\{a^1,b^1,c^1,d^1,e^1,f^1,a^2,b^2,c^2,d^2,e^2,f^2 \} \]
isomorphic to the complement of the
collinearity graph of a grid as in Fig. 3. Clearly, $D$ admits a
bipartition into two classes of size $6$ (corresponding to the two
long lines of the grid) as well as a $6$-partition into pairs,
corresponding to the six short lines. Two points in the same pair
are said to be {\it opposite.} Clearly, any four (respectively, three) pairs of
opposite vertices of a dodecade form a cube (a hexagon). We say that a dodecade  is {\it centric} if all the projective lines through opposite vertices are concurrent at the same point called the {\it center of the dodecade}.

\medskip

Let $\varphi=\varphi_0\circ \hat{\alpha}$ where $\varphi_0\colon \widehat{\X} \rightarrow \X$ is the canonical
$2$-fold covering of $\X$ as in Theorem~\ref{Main thm 1} and $\hat{\alpha}\in Aut(\widehat{\X}).$

For the subgraph $\Gamma_A$ of $\Gamma$ induced on a subset $A$
of the vertex set of $\X$, the subgraph $\widehat{\Gamma}_{A}$ of
$\widehat{\Gamma}$ induced on the subset $\varphi^{-1}(A)$ of
the vertex set of $\widehat{\Gamma}$ is called the {\it preimage} of
$\Gamma_A$ by $\varphi$.

Given a path $(X_0,X_1,\dots,X_n)$ in $\Gamma$ and a vertex $x_0\in \varphi^{-1}(X_0)$, there exists a unique path $(x_0,x_1,\dots, x_n)$ in $\widehat{\Gamma}$ such that $\varphi(x_i)=X_i$ for every $i=0,1,\dots, n.$ In graph theory and topology the path $(x_0,x_1,\dots, x_n)$ is called the {\it lifting} of $(X_0,X_1,\dots,X_n)$ {\it at} $x_0$ through $\varphi.$ By a harmless terminological abuse, we adapt this terminology to the induced subgraph $\Gamma_A$ of $\Gamma,$ calling {\it liftings} of $\Gamma_A$ the connected components of $\widehat{\Gamma}_A.$


The following theorem gives a characterization of the cliques of
$\Gamma$ in terms of their liftings.
Recall that $n_0$ is the nucleus of $Q_0.$

\begin{theorem}\label{Main thm 3}  Let $\varphi=\varphi_0\circ \hat{\alpha}$  with $\varphi_0$ and $\hat{\alpha}$ as above. \begin{itemize}
\item[(a)] The covering $\varphi$ induces a bijection between the set of  all centric hexagons of $\widehat{\Gamma}$ with center $n_0$ and the set of all non-linear $3$-cliques in $\Gamma.$
\item[(b)] The covering $\varphi$ induces a bijection between the set of  all centric cubes of $\widehat{\Gamma}$ with center $n_0$ and the set of all non-linear $4$-cliques in $\Gamma.$
\item[(c)] The covering $\varphi$ induces a bijection between the set of all centric dodecades of $\widehat{\Gamma}$ with center $n_0$ and the set of all non-linear $6$-cliques in $\Gamma.$

\end{itemize}
\end{theorem}

Let $S$ be a subset of $\PG(V).$ According to~\cite{DBP}, there exists a unique
subfield $\F_S$ of $\F_q$ such that $S$ is contained in a projective
subgeometry of $\PG(V)$ defined over $\F_S$ and $\F_S$ is minimal
with respect to this property. Moreover, the
family of projective subgeometries of $\PG(V)$ defined over $\F_S$
and containing $S$ admits a smallest member which we shall denote by
$\langle S\rangle_{\F_S}.$ We call $\langle S\rangle_{\F_S}$ the
{\it $\F_S$-span of $S$}.

If $S$ is a subset of  $Q,$ we define the {\it subgeometry} of $Q$ {\it induced} on $S$ to be the pair $Q(S):=(S,
\mathcal{L}_S)$ where

$\mathcal{L}_S:=\{l\cap S\colon l {\rm \,\,is\,\,a
\,\,line\,\,of\,\,Q\,\,such\,\,that\,\, } |l\cap S|\geq 2\}.$

As recalled a few lines above, there exists a unique smallest projective subgeometry $\langle  S\rangle_{\F_S}$ of $\PG(V)$ containing $S$ and defined over  a subfield $\F_S$ of $\F_q.$ We say that $Q(S)$ is an {\it $\F_2$-subgeometry} of $Q$ if $\F_S=\F_2.$

 The next theorem, combined with Theorem~\ref{Main thm 3},  gives a characterization of non-linear $3$-, $4$- and $6$-cliques in terms of $\F_2$-subgeometries of $Q.$

\begin{theorem} \label{Main thm 4} All the following hold.
\begin{itemize}
\item[1)] A centric hexagon of $\widehat{\X}$ with center $n_0$ is the subgraph of $\widehat{\Gamma}$ induced on the set-complement
$Q(S)\setminus H_0 $ of  $H_0=n_0^{\perp}$  in  an
$\F_2$-subgeometry $Q(S)\cong Q^+(3,2)$ of $Q$ such
that $n_0\in \langle S\rangle_{\F_2}.$
\item[2)] A centric cube of $\widehat{\X}$ with
center $n_0$ is the subgraph of $\widehat{\Gamma}$ induced on the set-complement $Q(S)\setminus H_0$  of
$H_0=n_0^{\perp}$  in an $\F_2$-subgeometry $Q(S)\cong Q(4,2)$ of $Q$
such that $n_0\in \langle S\rangle_{\F_2}.$
\item[3)] A centric
dodecade of $\widehat{\X}$ with center $n_0$ is the subgraph of $\widehat{\Gamma}$ induced on the set-complement $Q(S)\setminus
H_0$ of $H_0=n_0^{\perp}$  in a $\F_2$-subgeometry $Q(S)\cong
Q^-(5,2)$ of $Q$ such that $n_0\in \langle S\rangle_{\F_2}.$
\end{itemize}
\end{theorem}

The relevance of the parity of $n$ for the existence of maximal non-linear $6$-cliques (Theorem~\ref{Main thm 2}) is justified by Theorem~\ref{Main thm 3} and Theorem~\ref{Main thm 4}. Indeed, every non-linear  $6$-clique arises from the complement of $H_0$  in a $(2,4)$-subquadrangle $Q^-(5,2)$ of $Q$ (claim 3 of Theorem~\ref{Main thm 4}), but it is well known that $Q^-(5,2)$ is a subgeometry of $Q^-(5,2^n)$ if and only if  $n$ is odd.

As far as non-linear 4-cliques are concerned, each non-linear 4-clique arises from the complement of a tangent hyperplane in the (2,2)-subquadrangle $Q(4,2)$ of $Q^-(5,2^n).$ Observe that $Q(4,2)$ is a subgeometry of $Q^-(5,2^n)$ for any $n$. So, non-linear $4$-cliques exist regardless the parity of $n.$

\section{The semi-partial geometry $\X$}\label{sec 3}

We will keep the same terminology and notation as in Section~\ref{Introduction}.

\subsection{The quotient geometry $\widehat{\X}/R$}\label{quotient geometry}
Let $R$ be the equivalence relation defined on the pointset
$\widehat{P}$ of the affine generalized quadrangle $\widehat{\X}$ as follows: for
$x,y \in \widehat{P}$, put $xR y$ if and only if either $x=y$ or
$x^{\perp}\cap H_0 =y^{\perp}\cap H_0 (=x^{\perp}\cap y^{\perp})$ (see~\cite[8.4]{Pas-Diagram geometry}), where $\perp$ is the orthogonality relation defined by $Q$, as in Section~\ref{Introduction}. Note that distinct points of $\widehat{\X}$
correspond in $R$ if and only if they have distance $3$
in the collinearity graph of $\widehat{\X}.$

For $x\in \widehat{P}$ and $\hat{l}\in \widehat{L}$, let $[x]=\{y\in \widehat{P}\colon xR y\}$ and  $[\hat{l}]:=\{[x]\colon x\in \hat{l}\}.$  Given any two distinct lines  $\hat{l}$ and $\hat{m}$ of $\widehat{\X},$ we put $\hat{l}R \hat{m}$ if and only if $[\hat{l}]=[\hat{m}].$

We define the {\it quotient geometry}  $\widehat{\X}/R$ as the
point-line geometry having  $\{[x]\colon x\in \widehat{P}\}$ as the pointset and $\{[\hat{l}]\colon \hat{l}\in \widehat{L}\}$ as the lineset.  A point $[x]$
is incident with a line $[\hat{l}]$ if $[x]\in [\hat{l}].$


We will now give a different equivalent description of the quotient geometry $\widehat{\X}/R$. We recall some well-known facts about the action of the polarity defined by $\perp
$ on the subspaces of $\PG(V)$ (see e.g.~\cite{HT91}).
\begin{itemize}
\item[$\bullet$] If $x$ is a point of $Q,$ then $x^{\perp}$ is a hyperplane tangent to $Q$ at $x$.
\item[$\bullet$] If $S$ is a $3$-dimensional subspace of $\PG(V)$ such that $S\cap Q$ is an elliptic quadric, then $S^{\perp}$ is a line secant to $Q$.
\item[$\bullet$] If $S$ is a $3$-dimensional subspace of $\PG(V)$ such that $S\cap Q$ is a hyperbolic quadric, then $S^{\perp}$ is a line disjoint from $Q$.
\item[$\bullet$] If $\pi_p$ is a projective plane of $\PG(V)$ such that $\pi_p\cap Q=\{p\},$ then $\pi_p^{\perp}$ is a projective plane intersecting $Q$ in two lines meeting at $p$. Moreover $\pi_p\cap \pi_p^\perp =\{p\}$ and $\langle \pi_p,\pi_p^\perp\rangle=p^\perp.$
\end{itemize}

 Let $\nu$ be the unique non-trivial elation of $\PG(V)$ stabilizing $Q$ and having the nucleus $n_0$ of $Q_0$ as the center and the hyperplane $H_0=n_0^{\perp}$ as the axis. Clearly, $\nu$ is an involution.

The relation $R$ can be equivalently defined by means of the
collineation $\nu.$ More precisely, the orbits of $\nu$ on the
pointset of $\widehat{\X}$ are precisely the equivalence
classes of $R$ on $\widehat{P}$. Hence, for any point $x$ of $\widehat{P},$ the
equivalence class $[x]$ is precisely the set $\{x, \nu(x)\}$
obtained as the intersection $\langle x, n_0\rangle \cap Q$, where
$\langle x, n_0\rangle $ is the line (secant to $Q$) of $\PG(V)$ through $x$ and $n_0$.
Accordingly, a line $[\hat{l}]$ of $\widehat{\X}/R$ corresponds to the set
$(l\cup  \nu(l)) \setminus (l\cap \nu(l))$ obtained as the
intersection $\langle l, \nu(l)\rangle\cap Q$ minus the point $l\cap \nu(l),$ where $l=\hat{l}\cup \{l^\infty\}$ (see~(\ref{l infinity})). Note also that $l\cap \nu(l)=l^\infty$ and $\langle l,\nu(l)\rangle=\langle l, n_0\rangle.$

\subsection{Proof of Theorem~\ref{Main thm 1}}\label{Proof of Main Theorem 1}
According to the definition of $t$-fold covering given in Section~\ref{main results-sec 1} and the facts recalled in Section~\ref{quotient geometry}, it is easy to see that the projection $\pi\colon \widehat{\X}\rightarrow \widehat{\X}/R$ mapping any point $x$ of $\widehat{\X}$ to the point $[x]$ of $\widehat{\X}/R$
and any line $l$ of $\widehat{\X}$ to the line $[l]$ of $\widehat{\X}/R$ is a 2-fold covering.  Hence
\begin{lemma}\label{theorem-quotient-1}
$\widehat{\X}$ is a $2$-fold cover of $\widehat{\X}/R.$
\end{lemma}

\begin{lemma}\label{theorem-quotient-sec 2}
The geometries $\widehat{\X}/R$ and $\X$ are isomorphic.
\end{lemma}
\pr
Let $\psi_0\colon \X\rightarrow \widehat{\X}/R$ be the map defined on the points and the lines of $\X$ as follows:
for $X\in \mathcal{E}$, the perp $X^\perp$ of $X$ in $\PG(5,q)$ is a line of $\PG(5,q)$ through $n_0,$ secant for $Q.$ We set $\psi_0(X)=X^\perp\cap Q.$ For a line $\rho_p$ of $\X$, let $\pi_p$ be the unique tangent plane to all $X\in \rho_p$ (see Proposition~\ref{prop1}). Then ${\pi_p}^{{\perp}}$ is a plane of $\PG(5,q)$ intersecting $Q$ in two lines $l,l'$ such that $l\cap Q_0=l'\cap Q_0=l\cap l'=p.$ Then  $\psi_0(\rho_p)=[\hat{l}]=[\hat{l}'].$
It is easily seen that $\psi_0$ is an isomorphism of geometries. \eop
\medskip

Theorem~\ref{Main thm 1} follows from Lemma~\ref{theorem-quotient-1} and Lemma~\ref{theorem-quotient-sec 2}, in particular $\varphi_0:=\psi_0^{-1}\circ \pi $ where $\varphi_0$ is the canonical covering of $\X$ as defined in Theorem~\ref{Main thm 1}.

 \section{The collinearity graph of $\mathcal{X}$}~\label{collinearity graph}
In this section we will investigate properties of the collinearity
graph $\Gamma$ of $\X.$
Recall that the vertices of $\Gamma$ are the elliptic ovoids of $Q_0$ and two vertices $X_1$ and $X_2$
are adjacent if and only if they are tangent, i.e. $|X_1\cap X_2|=1.$

\begin{prop}
If $q>2$ then the collinearity graph $\Gamma$ of $\X$ has diameter $2$.
\end{prop}
\pr
The collinearity graph $\widehat{\Gamma}$ of $\widehat{\X}$ has diameter $3$ and the classes of $R$ are the pairs of points of $\widehat{\X}$ at mutual distance $3$ (see~\cite[8.4]{Pas-Diagram geometry}). Hence $\widehat{\Gamma}/R$ has diameter at most $2$. In fact, it is easy to see that it has diameter exactly $2$ except when $q=2$. By Lemma~\ref{theorem-quotient-sec 2} the graph $\Gamma$ has diameter $2$.
\eop

\medskip

If $q=2$, then $\Gamma$ is a complete graph on $6$ vertices. In this case $\widehat{\Gamma}$ is the set complement of a $(2\times 6)$-grid and the classes of $R$ are the short lines of the grid. We assume that $q\geq 4.$

In~\cite{HM-1981} it is proved that $\Gamma$ is a strongly regular graph with parameters

\noindent $v=q^2(q^2-1)/2,$ $k=(q-1)(q^2+1),$ $\lambda=(q-1)(q+2),$ and $\mu=2q(q-1).$

Relying on properties of the canonical covering $\varphi_0$ we can immediately obtain the parameters $(v,k,\lambda,\mu)$ of $\Gamma.$ Indeed,
\begin{itemize}
\item[$v:$] The number of vertices of $\Gamma$ is $|\widehat{P}|/2.$
\item[$k:$] The number of vertices adjacent to a given vertex $X$ of $\Gamma$ is

$|\{{\rm lines\,\,of\,\,} \widehat{\X}\,\,{\rm through\,\,} x\}| \times $

 $\times |\{{\rm points\,\,different \,\,from \,\,} x {\rm \,\,on\,\,any\,\,line\,\,of\,\,}\widehat{\X}{\rm \,\,through}\,\,x \}|,$

 for a point $x\in X^\perp \cap \widehat{P},$ no matter which.

\item[$\lambda:$] Let $X$ and $Y$ be two adjacent vertices of $\Gamma$, i.e. $X$ and $Y$ are two elliptic ovoids of $Q_0$ tangent at the point $p\in Q_0.$ Hence they define a line $\rho_p$ of $\X$. Then $\varphi_0^{-1}(\rho_p)=\pi_p^\perp\cap Q=\{l^1, l^2\}\setminus \{p\}$ where $l^1, l^2$ are lines of $Q$ through $p.$ Hence $\varphi_0^{-1}(X)=\{x^1, x^2\}\subset \varphi_0^{-1}(\rho_p)$ and $\varphi_0^{-1}(Y)=\{y^1, y^2\}\subset \varphi_0^{-1}(\rho_p).$ Suppose $x^1$ and  $y^1$ are on $\hat{l}^1$ and $x^2$ and  $y^2$ are on $\hat{l}^2.$ The number of vertices of $\Gamma$ adjacent to both $X$ and $Y$ is is the same as the number of points $z\in \widehat{P}$ on $\hat{l}^1\setminus \{x^1, y^1\}$ (equivalently, the number of points on $\hat{l}^2\setminus \{x^2, y^2\}$). These points bijectively correspond to the ovoids of $\rho_p$ different from $X$ and $Y$. The number of such ovoids is $q-2$. We must add the number of points $z\in \widehat{\X}$ at distance $1$ from both $x^1$ and $y^2$ (equivalently, at distance $1$ from both $x^2$ and $y^1$). There are $q^2$ such points. Indeed, on any line of $\widehat{\X}$ through $x^1$ different from $\hat{l}^1$ there is a unique point of $\widehat{\X}$ collinear to $y^2.$  Hence $\lambda=q^2+ q-2$.
\item[$\mu:$] Let $X$ and $Y$ be two non-adjacent vertices of $\Gamma$, i.e. $X$ and $Y$ are two elliptic ovoids of $Q_0$ intersecting in a conic. Suppose $\varphi_0^{-1}(X)=\{x^1, x^2\}$ and $\varphi_0^{-1}(Y)=\{y^1, y^2\}.$ The number of vertices of $\Gamma$ adjacent to both $X$ and $Y$ is the same as the number of points $z\in Q\setminus Q_0$  not collinear with any point of $X\cap Y$, which are collinear to $x^1$ and $y^1$ or to $x^2$ and $y^1.$  Let $\hat{l}$ be a line of $\widehat{\X}$  through $y^1$ with $l^\infty \in Y\setminus X.$ There are $q^2-q$ such lines. For each of these $q^2-q$ lines there are exactly two distinct points adjacent to only one of $x^1$ or $x^2.$ Hence $\mu=2(q^2-q).$

    Note that there are no elliptic ovoids $Z$ tangent at $X$ and $Y$ at a point in $X\cap Y$ because otherwise $X$ and $Y$ would be on the same line of $\X$, which is not possible because $X$ and $Y$ are not adjacent in $\Gamma.$
\end{itemize}

In the above calculation of the parameters $\lambda$ and $\mu$ we have also proved two interesting properties of the geometry $\X$ which we will state as propositions:

\begin{prop}
Let $X$ and $Y$ be two tangent elliptic ovoids of $Q_0.$  Take a point $x\in X_1\setminus X_2.$ Then there exists a unique elliptic ovoid of $Q_0$ through $x$ which is tangent to both $X_1$ and $X_2.$
\end{prop}

\begin{prop}
Let $X_1$ and $X_2$ be two elliptic ovoids of $Q_0$ intersecting in a conic.  If $x\in X_1\setminus X_2$ then there exist exactly two distinct elliptic ovoids of $Q_0$ through $x$ which are both tangent to $X_1$ and $X_2.$ If $x\in X_1\cap X_2$ then there exist no elliptic ovoids of $Q_0$ through $x$ which are both tangent to $X_1$ and $X_2.$
\end{prop}



\subsection{Cliques of $\Gamma$ and their liftings}
In this section we classify cliques in $\Gamma.$ A clique in
$\Gamma$ consisting of members of $\X$ sharing the same point
$p$ is called a {\it linear clique (based at $p$)}. A rosette of
$\X$ based at $p$, being a maximal pencil of elliptic ovoids mutually tangent at $p$, is a maximal linear clique (this is part of claim 1 and 2 of Theorem~\ref{Main thm 2}).

Let $\widehat{\Gamma}$ denote the collinearity graph of the
affine quadrangle $\widehat{\X}$, as in Section~\ref{Introduction}. We recall that $\widehat{\Gamma}$ is regular with
diameter $3,$ it contains $4$-circuits and every clique of $\widehat{\Gamma}$ is contained in a line of $\widehat{\X}$. Let $\varphi_0\colon
\widehat{\X}\rightarrow \X$ be the canonical covering as in Theorem~\ref{Main thm 1}. The fibers of $\varphi_0$ are the pairs of points of $\widehat{\X}$ at mutual distance $3$ (see Section~\ref{quotient geometry}). We will adopt the following notation. If $X$ is a vertex
of $\Gamma$, then $\varphi^{-1}(X)=\{x^1,x^2\} (=X^\perp \cap Q).$



 \subsection{Preimages of edges and linear cliques}
The preimage of an edge in $\Gamma$ is the union of two disjoint edges in $\widehat{\Gamma}$ which are the liftings of that edge. Given a line $\rho_p$ of $\X$, we have $\varphi_0^{-1}(\rho_p)=\hat{l}\cup\hat{m}$ for two lines $\hat{l},$ $\hat{m}$ of $\widehat{\X}$ such that $l^\infty=m^\infty=p.$ The preimage of the subgraph $\Gamma_{\rho_p}$ induced by $\Gamma$ on $\rho_p$ is the subgraph induced by $\widehat{\Gamma}$ on $\hat{l}\cup\hat{m}.$ This graph is the disjoint union of two complete subgraphs, namely $\hat{l}$ and $\hat{m}$.

\begin{prop}\label{cliques-1}
Let $\mathcal{C}$ be an $i$-clique of $\Gamma$, $i\geq 4$, containing a linear $3$-clique. Then $\mathcal{C}$ is a linear clique.
\end{prop}
\pr
Let $A,B,C, X$ be four distinct vertices of $\C,$ with $A,B,C$ forming a linear subclique and let $\rho$ be the line of $\X$ on $\{A,B,C\}.$

By way of contradiction suppose that $X\not\in \rho.$

Let us lift the ordered closed configuration $(X,A,B,C,X)$ of $\Gamma$ starting at $x^1$ recalling that $A,B,C$ are on the same line $\rho$ of $\X.$
Since $\varphi_0$ is a covering, the point $x^1$ is adjacent in $\widehat{\Gamma}$ with exactly one preimage of $A$, say $a^1$.
Similarly, $ a^1$ is adjacent with exactly one preimage $ b^1$ of $B$ and
 $b^1$ is adjacent with exactly one preimage $c^1$ of $C$.
Again, $c^1$ is adjacent with exactly one preimage $x^2$ of $X$.
Moreover, since $\varphi_0$ is a covering, $\varphi_0^{-1}(\rho)$ contains exactly one line $\hat{l}$ of $\widehat{\X}$ incident with all of $a^1,b^1, c^1.$ On the other hand, $x^1, x^2\not\in \hat{l}$ as $X\not\in \rho.$

Since there are no triangles in $\widehat{\X}$, $x^1\not= x^2.$
 The path $(x^1, a^1,b^1, c^1,x^2)$ is the lifting of $(X,A,B,C,X)$ at $x^1.$

 Lift now the closed ordered path $(X,A,B,X)$ of $\Gamma$ starting at $x^1.$
Proceeding in the same way as above, we get a path $(x^1, a^1,b^1,x^3)$ where the vertex $x^3$ must be different from $x^1$ because there are no triangles in  $\widehat{\X}.$
Similarly, $(x^2, c^1,b^1,x^3)$ is the lifting of $(X,C,B,X)$ at $x^2.$ Hence $x^3\not=x^2.$ (see Fig.4).




\begin{picture}(310,108)(0,0)
 \put(60,90){\line(1,0){150}}
\put(220,88){$\hat{l}$}
\put(79,87){$\bullet$}
\put(79,94){$a^1$}
\put(130,87){$\bullet$}
\put(130,94){$b^1$}

\put(180,87){$\bullet$}
\put(182,94){$c^1$}

\put(81,40){\line(0,1){49}}
\put(132,40){\line(0,1){49}}
\put(182,40){\line(0,1){49}}

\put(79,38){$\bullet$}
\put(79,30){$x^1$}

\put(130,38){$\bullet$}
\put(130,30){$x^3$}

\put(180,38){$\bullet$}
\put(180,30){$x^2$}
\put(126,10){Figure 4}

\end{picture}

It follows that $\varphi_0^{-1}(X)$ contains at least three distinct points, namely  $x^1,$ $x^2,$ $x^3.$ This is a contradiction because $\varphi_0$ is a $2$-fold covering.\eop

\subsection{Liftings of non-linear  cliques}
We establish the following notation for projective points:
$a=[v]=[v_1,\dots, v_n]$ where $v=(v_1,\dots, v_n)$ is a given vector representing the projective point $a.$

Given two projective points $a$ and $b$, symbols as $ta+sb$ with $t$ and $s$ scalars are nonsense. However, when particular representatives $v$ and $w$ have been chosen for  $a$ and $b$, we will take the liberty of writing $ta+sb$ for $[tv+sw]$. Similarly, given a quadratic form $f$, we write $f(a)$ for $f(v).$ Of course, these are notational abuses but they are harmless as far as it is clear which is the representative vector chosen for a given point. In particular, when doing so, we forbid ourselves any rescaling of the given vectors.

In view of Proposition~\ref{cliques-1}, in the following sections we
only consider non-linear cliques of $\Gamma.$

Suppose that the $(q,q^2)$-generalized quadrangle $Q$ is represented, with respect to a given basis,
by the  equation $f(x)=0$ where
\begin{equation}\label{equation Q1}
f((x_1,x_2,x_3,x_4,x_5,x_6)):=x_1x_2+x_3x_4+x_5^2+x_5x_6+\lambda x_6^2
\end{equation}

\noindent with $\lambda$ a given element of $\F_q$ with $Tr(\lambda)=1.$ Recall that the trace function $Tr\colon \F_q\rightarrow \F_2,$ $Tr(x)=\sum_{i=0}^{n-1}x^{2^i},$ ($q=2^n$), is  an $\F_2$-linear surjective homomorphism and a polynomial $x^2+xy+\lambda y^2$ is irreducible in $\F_q$ if and only if $Tr(\lambda)=1.$

We may assume that the nucleus $n_0$ of $Q_0$ is represented by the vector $(0,0,0,0,0,1).$ Accordingly, $H=n_0^\perp$ is the hyperplane of equation $x_6=0.$

A set of four distinct points $\{a,b,c,d\}$ of $Q$ is a
{\it quadrangle} of $Q$ if it is the set of points of a $(1,1)$-subquadrangle of $Q.$ Hence
\begin{equation}
q(a)=q(b)=q(c)=q(d)=0;\,\,\, a\perp b\perp c\perp d\perp a;\,\,\, a\not\perp
c; \,\,\, b\not\perp d.
\end{equation}
We first need a rather technical lemma.
 \begin{lemma}\label{lemma-algebraic condition}
Let $\lambda, \mu$ be two elements of $\F_q.$ Suppose  $Tr(\lambda)=1$ and $\mu\not=1$.
Then the equation $x^2+xy+\lambda y^2 +\mu=1,$ in the unknowns $x,y,$ admits $q+1$ distinct solutions in $\F_q^2$.
\end{lemma}
\pr
If $y=0$ then the pair  $(x,y)=(\mu^{2^{n-1}}+1,0)$ is the unique solution of $x^2+xy+\lambda y^2 +\mu=1.$
If $y\not=0$ then the equation $x^2+xy+\lambda y^2 +\mu=1$ is soluble if and only if the equation $t^2+t+\lambda +(\mu+1)/y^2 =0$ in the unknown $t:=x/y$ is soluble. This latter requirement is equivalent to $Tr(\lambda +(\mu+1)/y^2)=0.$  Since
\[Tr(\lambda +(\mu+1)/y^2)=Tr(\lambda) +Tr((\mu+1)/y^2)=1+Tr((\mu+1)/y^2),\]
the equation $t^2+t+(\lambda +(\mu+1)/y^2 )=0$ admits solutions if and only if $Tr((\mu+1)/y^2)=1.$

As, the trace function $Tr$  is a surjective homomorphism of $\F_2$-vector spaces there are $q/2$ distinct values of $y$ such that $Tr((\mu+1)/y^2)=1$ (see e.g.~\cite{LN}). For each of these $q/2$ values of $y$, there are two distinct solutions of $t^2+t+\lambda +(\mu+1)/y^2 =0.$ So, if $y\not=0$, there are $2\cdot \frac{q}{2}$ distinct pairs $(x,y)$ which satisfy $x^2+xy+\lambda y^2 +\mu=1.$
\eop

\subsubsection{Non linear 3-cliques and 4-cliques}\label{subsubsection 3-cliques and 4-cliques}
Hexagons, cubes and dodecades of $\widehat{\X}$ and their centers, if any, have been defined in Section~\ref{main results-sec 1}. Clearly, those definitions can be immediately generalized to $Q,$ by replacing $\widehat{\Gamma}$ with the collinearity graph of $Q.$ In particular, a hexagon (cube, dodecade) of $Q$ is {\it centric} if the projective lines joining pairs of opposite points  are concurrent in one point, called the {\it center} of the hexagon (cube, dodecade).
These more general notions will be used from times to times in the sequel.

\begin{theo}\label{cliques-3}
The preimage in $\widehat{\Gamma}$ of a non-linear $3$-clique of $\Gamma$ is a centric
hexagon of $\widehat{\X}$ with the nucleus $n_0$ of $Q_0$ as the center.
\end{theo}
\pr Let $\C=\{A,B, C\}$ be a non-linear $3$-clique of
$\Gamma.$ Let us lift the ordered closed configuration $(A,B,C,A)$ of $\C$  starting at $a^1\in \varphi_0^{-1}(A).$
We get a path $(a^1,b^1,c^1,a^2)$ of  $\widehat{\Gamma}$, where $a^1\not= a^2$ since in $\widehat{\X}$ there are no proper triangles. Consider now the lifting of  $(A,B,C,A)$ at $a^2.$ We get a path $(a^2,b^2,c^2,a^3)$ where $b^2\not= b^1$ and $c^2\not= c^1$ because there are no proper triangles  in $\widehat{\X}$ and $a^3$ is a preimage of $A$, different from $a^2.$ As $|\varphi_0^{-1}(A)|=2,$ necessarily  $a^3\not= a^1.$ So, by pasting $\{a^1,b^1,c^1,a^2\}$ with $\{a^2,b^2,c^2,a^1\},$ we obtain a proper hexagon $H.$ As $\varphi_0$ is a $2$-fold covering, $H=\varphi_0^{-1}(\C),$ namely $\C$ lifts to $H.$

By definition of the canonical covering $\varphi_0,$ each of the pairs $\{a^1,a^2\}$, $\{b^1,b^2\}$, $\{c^1,c^2\}$ is on a secant line to $Q$ through $n_0.$\eop
\medskip

Clearly, every centric hexagon of $\widehat{\Gamma}$ with center $n_0$ is mapped by $\varphi_0$ onto a non-linear $3$-clique of $\Gamma.$ This observation combined with  Theorem~\ref{cliques-3} proves claim a) of Theorem~\ref{Main thm 3}.
\begin{theo}\label{cliques-2.a}
The preimage of a non-linear $4$-clique of $\Gamma$ is a centric
cube of $\widehat{\X}$ with the nucleus $n_0$ of $Q_0$ as the center.
\end{theo}
\pr
Let $\C=\{A,B,C,D\}$ be a non-linear
$4$-clique of $\Gamma.$ Let $\varphi_0^{-1}(A)=\{a^1,a^2\},$ $\varphi_0^{-1}(B)=\{b^1,b^2\},$ $\varphi_0^{-1}(C)=\{c^1,c^2\},$ $\varphi_0^{-1}(D)=\{d^1,d^2\}.$

Let us lift the ordered closed path $(A,B,C,D,A)$ of $\C$ starting at $a^1.$
We get a path  $\C^1=(a^1,b^1,c^1,d^1, a^i)$ where $i=1$ or $i=2.$ Suppose $i=2$, namely $\C^1=(a^1,b^1,c^1,d^1,
a^2).$ Then consider the lifting of the ordered closed path
$(A,B,C,A)$ starting at $a^1.$ We get either $(a^1,b^1,c^1,a^1)$ or
$(a^1,b^1,c^1,a^2).$ In the first case we obtain a proper triangle in $\widehat{\X},$ which cannot be. In the second case $(a^1,b^1,c^1,a^2)$ is a triangle in $\widehat{\X},$ again a contradiction.

Hence $i=1$, namely $\C^1=(a^1,b^1,c^1,d^1, a^1).$

Keeping in mind that the liftings to $\widehat{\Gamma}$ of a vertex-subset of $\Gamma$ cannot contain any
proper triangle of $\widehat{\X},$ we have the following:

\noindent a) the lifting of $(A,B,C,A)$ starting at $a^1$ is $(a^1,b^1,c^1,a^2)$;

\noindent b) the lifting of $(A,C,D,A)$ starting at $a^1$ is $(a^1,c^2,d^2,a^2)$;

\noindent c) the lifting of $(B,C,A,B)$ starting at $b^1$ is $(b^1,c^1,a^2,b^2)$;

\noindent d) the lifting of $(B,D,A,B)$ starting at $b^1$ is $(b^1,d^2,a^2,b^2)$;

\noindent e) the lifting of $(D,B,C,D)$ starting at $d^1$ is $(d^1,b^2,c^2,d^2).$



We finally get the cube as below (compare  Fig. 2).

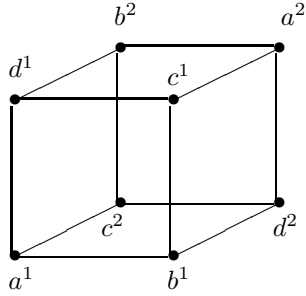
\begin{picture}(310,125)(0,0)
\put(0,18){$\bullet$}
\put(0,8){$a^1$}
\put(3,20){\line(1,0){57}}
\put(60,18){$\bullet$}
\put(60,8){$b^1$}
\put(1,20){\line(0,1){57}}
\put(61,20){\line(0,1){57}}
\put(0,88){$d^1$}
\put(60,85){$c^1$}
\put(3,80){\line(1,0){57}}
\put(0,77){$\bullet$}
\put(60,77){$\bullet$}
\put(2,20){\line(2,1){40}}
\put(62,20){\line(2,1){40}}
\put(2,79){\line(2,1){40}}
\put(62,79){\line(2,1){40}}
\put(40,38){$\bullet$}
\put(35,28){$c^2$}
\put(43,40){\line(1,0){57}}
\put(100,38){$\bullet$}
\put(100,28){$d^2$}
\put(41,40){\line(0,1){57}}
\put(101,40){\line(0,1){57}}
\put(40,108){$b^2$}
\put(103,108){$a^2$}
\put(43,100){\line(1,0){57}}
\put(40,97){$\bullet$}
\put(100,97){$\bullet$}
\end{picture}
\medskip

By definition of the covering $\varphi_0$, each of the pairs
$\{a^1, a^2 \}$, $\{b^1, b^2 \}$,
$\{c^1,  c^2 \}$, $\{d^1, d^2 \}$ is on a
secant line to $Q$ through $n_0.$ \eop
\medskip

Clearly, every centric cube of $\widehat{\Gamma}$ with center $n_0$ is mapped by $\varphi_0$ onto a non-linear $4$-clique of $\Gamma.$ This observation combined with  Theorem~\ref{cliques-3} proves claim b) of  Theorem~\ref{Main thm 3}.
\medskip

Let $\mathcal{F}$ be a quadrangle of $Q$. We denote by $I(\mathcal{F})$ the set of points $p\in \PG(V)\setminus Q$ such that $p^\perp \cap\mathcal{F}=\emptyset$ and there exists a cube in the affine quadrangle $Q\setminus p^\perp$ with center $p$ having $\mathcal{F}$ as one of its faces. Clearly, such a cube is uniquely determined by $\mathcal{F}$ and $p.$

\begin{lemma}\label{fundamental quadrangle}
Let $\bar{\mathcal{F}}=\{\bar{a}^1, \bar{b}^1, \bar{c}^1, \bar{d}^1\}$ be the quadrangle of $Q$ with
\medskip

\centerline{$\bar{a}^1=[1,0,0,0,0,0],\,\,\bar{b}^1= [0,0,1,0,0,0],\,\,\bar{c}^1=[0,1,0,0,0,0],\,\,\bar{d}^1= [0,0,0,1,0,0].$}



Then 
\begin{equation}
I(\bar{\mathcal{F}})=\{[u,\frac{1}{u},v,\frac{1}{v},r,s]\colon r^2+rs+ \lambda s^2=1,\, uv\not=0, u,v\in \F_q \}
\end{equation}
and
\begin{equation}\label{size I(F)}
|I(\bar{\mathcal{F}})|=(q-1)^2(q+1)=(q-1)(q^2-1).
\end{equation}
\end{lemma}
 \pr
Let $\bar{p}=[p_1,p_2,p_3,p_4,r,s]\in I(\bar{\mathcal{F}}).$
 Since $\bar{p}$  does not belong to $Q$, we have $f(\bar{p})\not=0$ (with $f(x)$ as in Equation~(\ref{equation Q1})). Without loss of generality we can assume $f(\bp)=1,$ i.e.
\begin{equation}\label{q(p)}
p_1p_2+p_3p_4+r^2+rs+\lambda s^2=1.
\end{equation}

 Since $\bar{p}^\perp \cap\bar{\mathcal{F}}=\emptyset$, $\bar{p}$ is not collinear in $Q$ with any element of $\bar{\mathcal{F}}.$  This is equivalent to $p_2\not=0,$ $p_4\not=0,$ $p_1\not=0,$ $p_3\not=0.$

Suppose that $\ba^2,$ $\bb^2,$ $\bc^2,$ $\bd^2$ are points of $Q$ so that $\bar{C}=\{\bar{a}^i,\bar{b}^i,\bar{c}^i,\bar{d}^i\}_{i=1,2}$ is a centric cube of $Q\setminus \bar{p}^\perp$ with center $\bp.$ Recall that the exponent $1$ and $2$ of the same symbol refer to opposite points, that the collinearities between points in $\bar{C}$ are as described in Fig. 2 and each of $\ba^2,$ $\bb^2,$ $\bc^2,$ $\bd^2$ are on a line joining $\bp$ to $\ba^1,$ $\bb^1,$ $\bc^1,$ $\bd^1$ respectively. Then, in terms of coordinates, we have

\[\ba^2=\ba^1+t_a \bp,\,\,\,\, \bb^2=\bb^1+t_b \bp,\,\,\,\, \bc^2=\bc^1+t_c \bp,\,\,\,\,\bd^2=\bd^1+t_d \bp\]

\noindent for suitable $t_a,t_b,t_c,t_d\in \F_q\setminus\{0\}.$ By hypothesis, plugging in the coordinates of $\ba^1,$ $\bb^1,$ $\bc^1,$ $\bd^1$, we have
\medskip

\noindent $\begin{array}{ll}
\ba^2=[1+t_ap_1,t_ap_2,t_ap_3, t_ap_4, t_a r, t_a s]; & \bb^2=[t_bp_1,t_bp_2,1+t_bp_3, t_bp_4, t_b r, t_b s];\\
 & \\
\bc^2=[t_cp_1,t_cp_2+1,t_cp_3, t_cp_4, t_c r, t_c s]; & \bd^2=[t_dp_1,t_dp_2,t_dp_3, 1+ t_dp_4, t_d r, t_d s].
\end{array}$
\medskip

Translating in terms of coordinates the collinearities in $Q$ for $\bar{C}$ to be a cube with center $\bp$, we have

\begin{equation}\label{conditions}
\left.\begin{array}{lll}
\!\!\!\!\!\!\!\!\!\!\!\!\ba^1\perp \bc^2:\, t_c p_2+1=0;& &  \bb^1\perp \bd^2:\, t_dp_4+1=0;\\
 & & \\
\!\!\!\!\!\!\!\!\!\!\!\!\bc^1\perp \ba^2:\, t_ap_1+1=0;& &  \bd^1\perp \bb^2:\, t_bp_3+1=0;\\
 & & \\
\!\!\!\!\!\!\!\!\!\!\!\!\ba^2\perp \bb^2:\, t_at_b p_2 p_1+t_a t_bp_3 p_4=0;& & \bb^2\perp \bc^2:\, t_bt_c p_1 p_2+t_b t_cp_3 p_4=0;\\
 & & \\
\!\!\!\!\!\!\!\!\!\!\!\!\bc^2\perp \bd^2:\, t_ct_d p_1 p_2+t_c t_dp_3 p_4=0;& & \bd^2\perp \ba^2:\, t_at_d p_1 p_2+t_a t_dp_3 p_4=0.\\
\end{array}\right\}
\end{equation}

The last four equations in~(\ref{conditions}) give $p_1p_2+p_3p_4=0.$
Hence, Equation~(\ref{q(p)}) becomes
\begin{equation} \label{conditions on r and s}
 r^2+rs+\lambda s^2=1.
 \end{equation}

Since $\ba^2,$ $\bb^2,$ $\bc^2,$ $\bd^2$ are points of $Q$, we have $q(\ba^2)=q(\bb^2)=q(\bc^2)=q(\bd^2)=0$, i.e. $p_2=t_a$,
$p_4=t_b$, $p_1=t_c$ and $p_3=t_d.$  The first four
equations in~(\ref{conditions}) give
$p_1p_2=1=p_3p_4.$

Put $u:=p_1$ and $v:=p_3$. Then $\bp=[u,1/u, v, 1/v,r,s]$ where $uv\not=0$ and $r,s$ satisfy Equation~(\ref{conditions on r and s}).

By Lemma~\ref{lemma-algebraic condition} with $\mu=0$, Equation~(\ref{conditions on r and s}) admits $q+1$ distinct solutions. The lemma is proved.\eop

\medskip

We will refer to the quadrangle $\bar{\mathcal{F}}$ as in Lemma~\ref{fundamental quadrangle} as the {\it fundamental quadrangle of $Q$}. By the proof of Lemma~\ref{fundamental quadrangle}, the coordinates of the points of the cube $C_{\bar{p}}$ having $\bar{\mathcal{F}}$ as a face and $\bar{p}\in I(\bar{\mathcal{F}})$ as the center are the following
\medskip
\begin{equation} \label{fundamental cube}
\left.\begin{array}{lll}
 \bar{a}^1=[1,0,0,0,0,0], & & \bar{b}^1=[0,0,1,0,0,0],\\
 \bar{c}^1=[0,1,0,0,0,0], & & \bar{d}^1=[0,0,0,1,0,0],\\
\ba^2=[0,\frac{1}{u^2} , \frac{v}{u} ,\frac{1}{uv} ,\frac{r}{u}, \frac{s}{u}]&  & \bb^2=[\frac{u}{v},\frac{1}{uv} , 0,\frac{1}{v^2} ,\frac{r}{v}, \frac{s}{v}];\\
\bc^2=[u^2,0,uv,\frac{u}{v},ur, us]&  & \bd^2=[uv,\frac{v}{u} , v^2,0, vr,vs]\\
\bar{p}=[u,1/u,v,1/v,r,s]. & & \\
\end{array}\right\}
\end{equation}
 where $uv\not=0$ and $r,s$ satisfy Equation~(\ref{conditions on r and s}).

\begin{co} Every quadrangle of $Q$ can be extended to $(q-1)^2(q+1)$ centric
cubes of $Q.$
 \end{co}
\pr Straightforward, by~(\ref{size I(F)}) of Lemma~\ref{fundamental quadrangle} and the transitivity of $SO^-(6,q)$ on the set of
quadrangles of $Q.$ 
\eop

\begin{co} $\Gamma$ admits non-linear $4$-cliques.
\end{co}
\pr
By Theorem~\ref{cliques-2.a},  if a non-linear 4-clique exists then it
can be uniquely lifted to a centric cube with center $n_0.$ Hence existence of non-linear 4-cliques is equivalent to the existence of centric cubes with center $n_0.$
By Lemma~\ref{fundamental quadrangle}, there exists centric cubes having $\bar{\mathcal{F}}$ as a face. Let $C_{\bp}$ one of these cubes  and let $\bp$ its center.

The group $SO^-(6,q)\leq Aut(Q)$ is transitive on the set of points of $PG(5,q)\setminus Q$. Then there exists $g\in SO^-(6,q)$  such that $g(\bp)=n_0$. The cube $g(C_{\bp})$ has the required properties.
\eop
\medskip

The following theorem is claim 3) of Theorem~\ref{Main thm 2}.
\begin{theo}\label{3-cliques extends to 4-cliques}
Every non-linear $3$-clique of $\Gamma$ can be extended to $q+1$ non-linear $4$-cliques.
\end{theo}
\pr By Theorem \ref{cliques-3} and \ref{cliques-2.a}, each non-linear $3$-clique can be extended to  $q+1$ non-linear 4-cliques if and only if each centric hexagon of $\widehat{\X}$ can  be extended to $q+1$ centric cubes of $\widehat{\X}$ having the same center as the hexagon.

Let $H=\{a^i, b^i,c^i\}_{i=1,2}$ be a centric hexagon of $Q$, with vertices marked as in Fig. 1.

Since the group $SO^-(6,q)$ is transitive on the set of paths of the collinearity graph of $Q$ of length $2$ we can suppose that $a^1$, $b^1$, $c^1$ are indeed the vertices $\ba^1$, $\bb^1$, $\bc^1$ as in the fundamental quadrangle $\bar{\mathcal{F}}$, i.e. $a^1$, $b^1$, $c^1$ are represented by the vectors $(1,0,0,0,0,0),$ $(0,0,1,0,0,0),$ and  $(0,1,0,0,0,0)$ respectively.

Let us determine the opposite points $a^2,$ $b^2$, $c^2.$ In order to do that, let $p=[p_1,p_2,p_3,p_4,p_5,p_6]$ be the center of $H$. Hence, $p$ is not collinear in $Q$ with any of $a^1,b^1, c^1.$ In terms of coordinates, this fact is equivalent to $p_1,p_2,p_3\not=0.$ Since $p$ is a point of $\PG(V)$ not in $Q$, $f((p_1,p_2,p_3,p_4,p_5,p_6))\not=0$. Without loss of generalities we can suppose $f((p_1,p_2,p_3,p_4,p_5,p_6))=1$, i.e. $p_1p_2+p_3p_4+p_5^2+p_5p_6+\lambda p_6^2=1.$ Then \[a^2=[\mu_a(1,0,0,0,0,0)+(p_1,p_2,p_3,p_4,p_5,p_6)]=[\mu_a+p_1,p_2,p_3,p_4,p_5,p_6];\]
\[b^2=[\mu_b(0,0,1,0,0,0)+(p_1,p_2,p_3,p_4,p_5,p_6)]=[p_1,p_2,\mu_b+p_3,p_4,p_5,p_6];\]
\[c^2=[\mu_c(0,1,0,0,0,0)+(p_1,p_2,p_3,p_4,p_5,p_6)]=[p_1,\mu_c+p_2,p_3,p_4,p_5,p_6].\]

Since $a^2$, $b^2$, $c^2$ are points of $Q$ and recalling that $f(p)=1$ by assumption, we get $\mu_a=1/p_2$; $\mu_b=1/p_4$; $\mu_c=1/p_1$, so
\begin{equation}\label{hexagon}\left.\begin{array}{l}
a^2=[1/p_2+p_1,p_2,p_3,p_4,p_5,p_6];\\
b^2=[p_1,p_2,1/p_4+p_3,p_4,p_5,p_6];\\
c^2=[p_1,1/p_1+p_2,p_3,p_4,p_5,p_6].
\end{array}\right\}
\end{equation}
Translating in terms of coordinates the collinearity relations in $H$ (see Fig.~1), we have
\[\frac{1}{p_2}+p_1=0\,\,\,\mbox{because}\,\,\, a^2\perp c^1\,\,\, {\rm and }\,\,\,\, \frac{1}{p_1}+p_2=0\,\,\,\mbox{because}\,\,\,  c^2\perp a^1.\]

Hence $p_2=1/p_1.$ Note that the relations $b^2\perp a^2$ and $b^1\perp c^1$ are automatically satisfied.

Hence we get for Equations~(\ref{hexagon}):
\begin{equation}\label{hexagon-1}\left.\begin{array}{l}
a^2=[0,1/p_1,p_3,p_4,p_5,p_6];\\
b^2=[p_1,1/p_1,1/p_4+p_3,p_4,p_5,p_6];\\
c^2=[p_1,0,p_3,p_4,p_5,p_6].
\end{array}\right\}
\end{equation}

and $p=[p_1, 1/p_1,p_3,p_4,p_5,p_6]$ with $p_3p_4+p_5^2+p_5p_6+\lambda p_6^2=0.$

 Let $d^1=[d_1,d_2,d_3,d_4,d_5,d_6]$ and $d^2=\mu d^1+p$ ($\mu\in \F_q$) be two points of $Q.$ If we require the set  $H\cup \{ d^1, d^2\}$ be a cube with center  $p$ we obtain the following conditions on the coordinates of $d^1$ and $d^2:$
$d_2=0$ because $d^1\perp a^1$ and $d_2=0$ because $d^1\perp c^1$. Also $d_4\not=0$ since $d^1\not\perp b^1.$ Hence we can assume $d_4=1.$ Therefore
$d^1=[0,0,d_3,1,d_5,d_6]$. Since $d^1\perp b^2$ we have
\begin{equation}\label{eq- ext-1}
d_3 p_4=\frac{1}{p_4}+p_3+d_5p_6+d_6p_5.
\end{equation}
Observe that $d^1\not\perp b^1.$ Indeed \[\alpha((p_1, 1/p_1,p_3,p_4,p_5,p_6),(0,0,d_3,1,d_5,d_6))=1/p_4\not=0,\]  where $\alpha$ is the bilinear form associated to $f.$
Since $d^2$ is a point of $Q$, we get
\begin{equation}\label{eq- ext-2}
d_3 =d_5^2+d_5d_6+\lambda d_6^2.
\end{equation}

Let us  turn to the opposite point $d^2= \mu d^1+ p$ of $d^1.$ The fact that $d^2$ is a point of $Q$ is equivalent to $1+\mu/p_4=0,$ hence $\mu=p_4.$
So, \[d^2=[p_1,1/p_1,p_4d_3+p_3,0,p_4d_5+p_5, p_4d_6+p_6].\] By (\ref{eq- ext-1}), we get
\[d^2=[p_1,1/p_1,p_5d_6+d_5p_6+1/p_4, 0,p_4d_5+p_5, p_4d_6+p_6].\]

Note that the relations $d^2\perp b^1$, $d^2\perp a^2$, $d^2\perp c^2$ are automatically satisfied.

By conditions~(\ref{eq- ext-1}) and~(\ref{eq- ext-2}) we now have
\[ \left(\frac{1}{p_4}\right)^2 + \frac{p_3}{p_4} + d_5\frac{p_6}{p_4} +d_6 \frac{p_5}{p_4} +d_5^2+d_5d_6+\lambda d_6^2=0\]
which is equivalent to

\[\left(d_5+\frac{p_5}{p_4}\right)^2+ \left(d_5+\frac{p_5}{p_4}\right)\left(d_6+\frac{p_6}{p_4}\right) + \lambda \left(d_6+ \frac{p_6}{p_4}\right)^2+\]
\[+\left(\frac{p_5}{p_4}\right)^2 + \frac{p_5 p_6}{p_4^2} + \lambda \left(\frac{p_6}{p_4}\right)^2 + \left(\frac{1}{p_4}\right)^2 + \left(\frac{p_3}{p_4}\right)=0.\]

Since \[\left(\frac{p_5}{p_4}\right)^2 + \frac{p_5 p_6}{p_4^2} + \lambda \left(\frac{p_6}{p_4}\right)^2 + \left(\frac{1}{p_4}\right)^2 + \left(\frac{p_3}{p_4}\right)=\]
\[\frac{1+p_3p_4 +p_5^2+p_5p_6+ \lambda p_6^2}{p_4^2}= \frac{f((p_1,p_2,p_3,p_4,p_5,p_6))}{p_4^2}=\frac{1}{p_4^2}\not=0\]
we can apply Lemma~\ref{lemma-algebraic condition} with $x:= d_5+\frac{p_5}{p_4}$, $y:= d_6+\frac{p_6}{p_4}$, $\mu= \frac{1}{p_4^2} +1$ to obtain $q+1$ different choices for $d^1.$ The theorem is proved.
\eop

\begin{prop}\label{cliques-2-1}The following hold.\begin{itemize}
\item[a)] The number of distinct non-linear $3$-cliques of $\Gamma $
 is  $N_3=\frac{(q^4-1)(q-1)q^4}{12}.$
\item[b)] The number of distinct non-linear $4$-cliques of  $\Gamma$ is  $N_4=\frac{(q^4-1)(q^2-1)q^4}{48}.$
\end{itemize}
\end{prop}
\pr
We first prove Claim b). By Theorem~\ref{cliques-2.a}, the lifting of a non-linear $4$-clique is a centric cube $C.$ Let us double count the set of pairs $\{(p, C)\}$ where $p$ is a point not in $Q$ and $C$ is a cube of $Q\setminus p^\perp$ with center $p.$ Hence

\centerline{$N_4\cdot$ ($\sharp$ of points not in  $Q$)$= N_4\cdot \frac{q^6-1-(q^3+1)(q^2-1)}{q-1}=N:=$ ($\sharp$ of centric cubes in  $Q$).}
\medskip

\noindent The total number $N$ of centric cubes in $Q$ can be computed as follows:

\medskip
$\begin{array}{l}
\rm {(\sharp \,\,of\,\, cubes\,\,  with\,\, a\,\, given\,\, quadrangle\,\,as\,\,a\,\,face)}  \cdot  \rm {(\sharp \,\,of\,\, quadrangles)}=\\
N  \cdot  \rm {(\sharp\,\, of\,\, quadrangles\,\, in\,\, a\,\, given\,\, cube).}
\end{array}$
\medskip

Counting the number of choices for the vertices of a quadrangle of $Q$, it is not difficult to see that there are $(q^3+1)(q^2+1)(q+1)q^6/8$ quadrangles in $Q.$ By Corollary \ref{fundamental quadrangle}, there are $(q-1)^2(q+1)$ centric cubes  through a given quadrangle. Since every cube has $6$ faces, we obtain $N=(q^3+1)(q^4-1)(q^2-1)q^6/48.$

It follows that the number of cubes of $Q$ with a given center is
\[N_4=N\cdot \frac{q-1}{(q^6-1)-(q^3+1)(q^2-1)}=\frac{(q^4-1)(q^2-1)q^4}{48}.\]

\noindent We now turn to claim a). By Theorem~\ref{cliques-3}, the lifting of a non-linear $3$-clique is a centric hexagon. Moreover, every non-linear $3$-clique  is contained in $q+1$ non-linear $4$-cliques, by Theorem~\ref{3-cliques extends to 4-cliques}. 
Hence
\medskip

$N_3\cdot$ ($\sharp$ of non-linear $4$-cliques on a non-linear $3$-clique)=\\
\indent $N_4\cdot$ ($\sharp$ of non-linear $3$-cliques in a non-linear $4$-clique).
\medskip

By Part b) of this Proposition, $N_4=\frac{(q^4-1)(q^2-1)q^4}{48}.$ Since there are four non-linear $3$-cliques contained in a $4$-clique, the proposition is proved.
\eop
\medskip

Note that claim a) of Proposition~\ref{cliques-2-1} also easily follows by the parameters $v, k, \lambda, \mu$ of $\Gamma$ as recalled at the beginning of Section~\ref{collinearity graph}.
 \subsubsection{Non linear 5-cliques and 6-cliques}
We have defined dodecades of $\widehat{\X}$ in Section~\ref{main results-sec 1}. Decades can be defined in a similar way. Explicitly, a {\it  decade} of $\widehat{\X}$ is a subgraph of
$\widehat{\Gamma}$ with $10$ vertices
\[T=\{a^1,b^1,c^1,d^1,e^1,a^2, b^2,c^2,d^2,e^2 \}\]
isomorphic to the complement of the collinearity
graph of a ($5\times 2$)-grid as follows:

\begin{picture}(310,100)(0,0)
\put(0,21){$\bullet$}
\put(0,8){$a^1$}
\put(3,24){\line(1,0){200}}

\put(53,21){$\bullet$}
\put(53,8){$b^1$}

\put(103,21){$\bullet$}
\put(103,8){$c^1$}

\put(153,21){$\bullet$}
\put(153,8){$d^1$}

\put(203,21){$\bullet$}
\put(203,8){$e^1$}

\put(1,23){\line(0,1){50}}
\put(54,23){\line(0,1){50}}
\put(104,23){\line(0,1){50}}
\put(154,23){\line(0,1){50}}
\put(204,23){\line(0,1){50}}
\put(0,71){$\bullet$}
\put(0,78){$a^2$}
\put(53,71){$\bullet$}

\put(53,78){$b^2$}

\put(103,71){$\bullet$}
\put(103,78){$c^2$}

\put(153,71){$\bullet$}
\put(153,78){$d^2$}

\put(203,71){$\bullet$}
\put(203,78){$e^2$}


\put(3,73){\line(1,0){200}}
\put(106,-6){Figure 5}
\end{picture}
\medskip

The graph $T$ admits a bipartition in two
classes of size $5$ (corresponding to the two long lines of the
grid) as well as a $5$-partition in pairs, corresponding to the five
short lines. Two points in the same pair are said to be {\it
opposite.} Clearly, any four (three) pairs of opposite vertices of a
decade form a cube (a hexagon). We say that a decade  is {\it
centric} if all the projective lines through opposite vertices are
concurrent at the same point called the {\it center of the decade}.

\begin{prop}\label{5-cliques}The following hold.
\begin{itemize}
\item[i)] Every non-linear $5$-clique of $\Gamma$ lifts through $\varphi_0$ to a centric
decade in $\widehat{\Gamma}$ with the nucleus $n_0$ of $Q_0$ as the center.
\item[ii)] Every non-linear $6$-clique of $\Gamma$ lifts through $\varphi_0$ to a a centric
dodecade of $\widehat{\Gamma}$ with the nucleus $n_0$ of $Q_0$ as the center.
\end{itemize}
\end{prop}

\pr
We will only prove claim ii). The proof of claim i) is similar.\\
Let $\C=\{a,b,c,d,e,f\}$ be a non-linear $6$-clique of $\Gamma.$
Consider the $4$-cliques $\{a,b,c,d\}$ and $ \{a,f,e,c\}$ of $\C$ and denote their liftings by $\{a^i,b^i,c^i,d^i\}_{i=1,2}$ and $\{a^i,f^i,e^i,c^i\}_{i=1,2}$, respectively. By Theorem~\ref{cliques-2.a},  the lifting of any non-linear $4$-subclique of $\C$ is a centric cube of $\widehat{\X}$ with center $n_0$. Hence we can lift $\C$ to a subgraph of $\widehat{\Gamma}$ which is the complement of the following graph

 \begin{picture}(310,100)(0,0)
\put(0,21){$\bullet$}
\put(0,8){$a^2$}
\put(3,24){\line(1,0){250}}

\put(53,21){$\bullet$}
\put(53,8){$b^1$}

\put(103,21){$\bullet$}
\put(103,8){$c^2$}

\put(153,21){$\bullet$}
\put(153,8){$d^1$}

\put(203,21){$\bullet$}
\put(203,8){$e^2$}

\put(253,21){$\bullet$}
\put(253,8){$f^1$}

\put(1,23){\line(0,1){50}}
\put(54,23){\line(0,1){50}}
\put(104,23){\line(0,1){50}}
\put(154,23){\line(0,1){50}}
\put(204,23){\line(0,1){50}}
\put(254,23){\line(0,1){50}}
\put(0,71){$\bullet$}
\put(0,78){$a^1$}
\put(53,71){$\bullet$}
\put(53,78){$b^2$}

\put(103,71){$\bullet$}
\put(103,78){$c^1$}

\put(153,71){$\bullet$}
\put(153,78){$d^2$}

\put(203,71){$\bullet$}
\put(203,78){$e^1$}

\put(253,71){$\bullet$}
\put(253,78){$f^2$}

\put(3,73){\line(1,0){250}}
\put(114,-6){Fig. 6}

\end{picture}
\medskip

Clearly, the complement of the graph in Fig. 6 is a dodecade. By definition of the canonical covering $\varphi_0$, each of the opposite pairs  $\{a^1, a^2\}$, $\{b^1, b^2\}$, $\{c^1, c^2\}$, $\{d^1, d^2\},$ $\{e^1, e^2\}$, $\{f^1, f^2\}$ is on a secant line to $Q$ through $n_0.$ \eop
\medskip

Every centric dodecade of $\widehat{\Gamma}$ with center $n_0$ is mapped by $\varphi_0$ onto a non-linear $6$-clique of $\Gamma.$ This observation combined with claim ii) of  Proposition~\ref{5-cliques} proves claim c) of  Theorem~\ref{Main thm 3}.
\medskip

The following theorem comprises  claims 1) and 2) of Theorem~\ref{Main thm 2}.

 \begin{theo}\label{cliques-4} The following hold.
\item[a)] The graph $\Gamma$ admits non-linear $5$-cliques if and only if $q=2^n$, $n$ odd.
\item[b)] Let $q=2^n$, $n$ odd. Then every $5$-clique is contained in a $6$-clique.
\item[c)] There exist no $k$-cliques with $k>6.$
\end{theo}
\pr Let $C:=\{a^i,b^i,c^i,d^i\}_{i=1,2}$
be the lifting of a non-linear $4$-clique of $\Gamma.$ By Theorem
\ref{cliques-2.a}, $C$ is a centric cube of $\widehat{\X}.$  By the transitivity of the action of $SO^-(6,q)$ on quadrangles of $Q$ and by the proof of Lemma~\ref{fundamental quadrangle} we can suppose
without loss of generalities to choose the cube $C_{\bp}$ as in (\ref{fundamental cube}) of Section~\ref{subsubsection 3-cliques and 4-cliques}. In particular, $C_{\bp}$ has center $\bp=[u,1/u, v, 1/v,r,s] $ with $uv\not=0$ and $r^2+rs+\lambda s^2=~1.$



We want to enlarge $C_{\bp}$ in such a way to get a lifting of a $k$-clique of $\Gamma$ with $k\geq 5.$

Let $\be^1=[e_{1,i}]_{i=1}^6$ and $\be^2=[e_{2,i}]_{i=1}^6$ be two points of $Q.$
The conditions on the coordinates of $\be^i$ ($i=1,2$) for the set $C_{\bp}\cup \{\be^i\}_{i=1}^2$ to be a decade of $Q\setminus \bp^\perp$ with center $\bp$ are the following.

\begin{equation}\label{conditions-5-clique}
\noindent \left.\begin{array}{lll}
\!\!\!\!\!\!\!\!\!\!\!\!\!\!\!\!\bb^1\perp \be^1:~~ e_{1,4}=0;&  & \bd^1\perp \be^1:~~ e_{1,3}=0;\\
\!\!\!\!\!\!\!\!\!\!\!\!\!\!\!\! & & \\
\!\!\!\!\!\!\!\!\!\!\!\!\!\!\!\!\ba^2\perp \be^1:~~ e_{1,1}=e_{1,5}us+e_{1,6}ur;&  & \bc^2\perp \be^1:~~ e_{1,2}=e_{1,5}s/u+e_{1,6}r/u;\\
\!\!\!\!\!\!\!\!\!\!\!\!\!\!\!\! & & \\
\!\!\!\!\!\!\!\!\!\!\!\!\!\!\!\!\ba^1\perp \be^2:~~ e_{2,2}=0;&  & \bc^1\perp \be^2:~~ e_{2,1}=0;\\
\!\!\!\!\!\!\!\!\!\!\!\!\!\!\!\! & & \\
\!\!\!\!\!\!\!\!\!\!\!\!\!\!\!\!\bb^2\perp \be^2:~~ e_{2,3}=e_{2,5}s v+e_{2,6}r v;&  & \bd^2\perp \be^2:~~ e_{2,4}=e_{2,5}s/v+e_{2,6}r/v .\\
\end{array}\right\}
\end{equation}


Since $\be^1$ and $\be^2$ are points of $Q$, we have $f(\be^1)=0$ and $f(\be^2)=0.$ Hence

\begin{equation}\label{cond-5-cliques-0}
e_{1,5}^2(s^2+1)+e_{1,5}e_{1,6}+e_{1,6}^2(r^2+\lambda)=0,
\end{equation}
\begin{equation}\label{cond-5-cliques}
e_{2,5}^2(s^2+1)+e_{2,5}e_{2,6}+e_{2,6}^2(r^2+\lambda)=0,
\end{equation}

The points $\be^1, \be^2$ and $\bp$ must be on the same line. Hence
\begin{equation}\label{cond-5-cliques-1}
e_{2,5}=e_{1,5}(1+rs)+e_{1,6}r^2\,\,\,  { \rm and}\,\,\,\ e_{2,6}=e_{1,5}s^2+ e_{1,6}(1+rs).
\end{equation}

By plugging $e_{2,5}$ and $e_{2,6}$ of Equation~(\ref{cond-5-cliques-1}) in Equation~(\ref{cond-5-cliques}) we have
 \[e_{1,5}^2((r^2+rs+\lambda s^2)s^2+1)+e_{1,5}e_{1,6}+e_{1,6}^2(r^2(r^2+rs+\lambda s^2)+\lambda)=0.\]
 Since $r^2+rs+\lambda s^2=1,$ the above equation gives us back Equation (\ref{cond-5-cliques-0}). So, Equation~(\ref{cond-5-cliques}) is already implicit in Equation~(\ref{cond-5-cliques-0}) hence we can disregard it.

The coordinates  $e_{1,i}$ of $\be^1$ for $1\leq i\leq 4$ and the coordinates $e_{2,i}$ of $\be^2$ for $1\leq i\leq 4$ depend on $e_{1,5}$ and $e_{1,6}.$ Hence we write  $\be^1_{(e_{1,5},e_{1,6})}$ for $\be^1$ and $\be^2_{(e_{1,5},e_{1,6})}$  for $\be^2$ so that to emphasize the dependance of $\be^1$ and $\be^2$ on $e_{1,5}$ and $e_{1,6}.$ 

\noindent Explicitly, we have
\begin{equation}\label{coord e1}
\!\!\!\!\!\!\!\!\!\!\!\!\!\!\!\!\!\be^1_{(e_{1,5},e_{1,6})}=(e_{1,5}us+e_{1,6}ur, e_{1,5} \frac{s}{u} +e_{1,6}\frac{r}{u}, 0,0,e_{1,5}, e_{1,6})
\end{equation}
and 
\begin{equation}\label{coord e2}
\!\!\!\!\!\!\!\!\!\!\!\!\!\!\!\!\!\be^2_{(e_{1,5},e_{1,6})}\!\!=\!\!(0,0,e_{1,5}sv+e_{1,6}rv, \frac{se_{1,5}}{v}+\frac{re_{1,6}}{v},e_{1,5}(1+rs)+e_{1,6}r^2\!\!, e_{1,5}s^2+e_{1,6}(1+rs)).
\end{equation}

Note that the possibility  to enlarge $C_{\bp}$ in order to get a lifting of a non-linear $k$-clique with $k\geq 5$ is equivalent to the existence of at least two points $\be^1$ and $\be^2$ as above. The number of choices for the pair $(\be^1,\be^2)$ determines the size of the $k$-clique. After having implemented all the collinearity relations among the points $C_{\bp}\cup \{\be^1, \be^2\}$ to get a centric decade, we got that the points $\be^1$ and $\be^2$ exist provided their coordinates satisfy Equations~(\ref{conditions on r and s}) and~(\ref{cond-5-cliques-0}). In the sequel we will show that Equations~(\ref{conditions on r and s}) and~(\ref{cond-5-cliques-0}) are equivalent to $Tr(1)=1.$ Then,  $Tr(1)=1$ if and only if $q$ is an odd power of $2$ (see e.g.~\cite{LN}).

More precisely, we will deal with the cases $s=1$ and $s\not=1$ separately showing that 

\begin{itemize}
\item[(i)] for $s=1$ and $s\not= 1$ the choices for $(\be^1,\be^2)$ allow exactly two extensions of $C_{\bp}$ to decades;
\item[(ii)] the union of the two decades in (i) is a centric dodecade.
\end{itemize}


Suppose firstly $s=1$. Then Equation~(\ref{cond-5-cliques-0}) is
\begin{equation}\label{eq s=1}
e_{1,5}e_{1,6}+e_{1,6}^2(r^2+\lambda)=0
\end{equation}

and Equation~(\ref{conditions on r and s}) becomes $\lambda=r^2+r+1.$

So, $Tr(\lambda)=Tr(r^2+r+1)=Tr(1)$. Since $Tr(\lambda)=1,$
\begin{equation}\label{condition-trace}
Tr(1)=1.
\end{equation}

Suppose now $s\not=1$. Then Equation~(\ref{cond-5-cliques-0}) is

\begin{equation}\label{cond-5-cliques-3}
e_{1,5}^2+ \frac{e_{1,5}e_{1,6}}{1+s^2}+\frac{(\lambda + r^2)(1+s^2)e_{1,6}^2}{(1+s^2)^2}=0.
\end{equation}

Equation~(\ref{cond-5-cliques-3}) admits solutions (in $e_{1,5}$ and $e_{1,6}$) if and only if

\begin{equation}\label{cond-5-cliques-4}
Tr(\lambda +\lambda s^2+r^2+r^2s^2)=0.
\end{equation}

By Equation~(\ref{conditions on r and s}) we have $\lambda s^2+r^2+r^2s^2=1+rs+r^2s^2$. Hence  Equation~(\ref{cond-5-cliques-4}) becomes $Tr(\lambda +\lambda s^2+r^2+r^2s^2)=Tr(\lambda + 1+rs+r^2s^2)=0$ which implies $Tr(\lambda) + Tr(1) + Tr(rs)+ Tr(r^2s^2)=0$, that is $Tr(1)=1$, again.  
\medskip

Points $\be^1$ and $\be^2$ exist if and only if Equation~(\ref{condition-trace}) is satisfied, which is equivalent to require  $q=2^n$ with $n$ odd.
So, non-linear $5$-cliques  exist if and only if $q=2^n$ with $n$ odd. Claim a) is proved.

Let $q=2^n$ and $n$ odd. Then Equations~(\ref{eq s=1}) and~(\ref{cond-5-cliques-3}) admit exactly two different solutions (in $e_{1,5}$ and $e_{1,6}$), up to proportionality: denote them by $(e_{1,5}^\star, e_{1,6}^\star)$ and $(e_{1,5}^\circ, e_{1,6}^\circ)$. According to~(\ref{coord e1}) and~(\ref{coord e2})  and keeping in mind that the exponents 1 and 2  refer to opposite points, put
\medskip

$\begin{array}{ll}
\be^{\star 1}:=\be^1_{(e_{1,5}^\star, e_{1,6}^\star)}, & \be^{\star 2}:=\be^2_{(e_{1,5}^\star, e_{1,6}^\star)}\\
  &  \\
\mbox{and} & \\
  & \\
\be^{\circ 1}:=\be^1_{(e_{1,5}^\circ, e_{1,6}^\circ)}, & \be^{\circ 2}:=\be^2_{(e_{1,5}^\circ, e_{1,6}^\circ)}.\\
\end{array}
$
\medskip

By construction, $C_{\bp}\cup\{\be^{\star 1},\be^{\star 2}\}$ and $C_{\bp}\cup\{\be^{\circ 1},\be^{\circ 2}\}$ are centric decades with center $\bp$ and $C_{\bp}\cup\{\be^{*1},\be^{*2}, \be^{\circ 1},\be^{\circ 2}\}$ is a centric dodecade with center $\bp.$ The Theorem is proved.
\eop
\medskip

The following is also implicit in the proof of Theorem~\ref{cliques-4} and corresponds to claims $4.i$ and $4.ii$ of Theorem~\ref{Main thm 2}.

\begin{co}\label{corollary-decades-dodecades}
Let $n$ be odd. Then every centric cube of $\widehat{\X}$ is contained in exactly two centric decades and just one centric dodecade of $\widehat{\X}.$ Every centric decade of $\widehat{\X}$ is contained in exactly one centric dodecade of $\widehat{\X}.$
\end{co}

\begin{prop}\label{number of 6-cliques}
Let $q=2^n,$ $n$ odd. The following hold.\begin{itemize}
\item[i)] The number of distinct non-linear $5$-cliques of $\Gamma$
 is  $N_5=\frac{(q^4-1)(q^2-1)q^4}{15\cdot8}.$\\
\item[ii)]The number of distinct non-linear $6$-cliques of $\Gamma$
is $N_6=\frac{(q^4-1)(q^2-1)q^4}{15\cdot48}.$
\end{itemize}
\end{prop}
\pr
By claim i) of Proposition~\ref{5-cliques}, the lifting of a non-linear $5$-clique is a centric decade with center $n_0.$ By Corollary~\ref{corollary-decades-dodecades}, the number of decades containing a given centric cube is two. Hence
$2\cdot N_4=N_5\cdot{5\choose 4}$, where $N_4=\frac{(q^4-1)(q^2-1)q^4}{48}$ is the number of centric cubes with center $n_0$ by claim b) of Proposition~\ref{cliques-2-1}. Therefore $N_5=\frac{(q^4-1)(q^2-1)q^4}{15\cdot8}.$

The lifting of a non-linear $6$-clique is a centric dodecade by claim ii) of Proposition~\ref{5-cliques}.
By Corollary~\ref{corollary-decades-dodecades}, the number of dodecades containing a given centric cube is one. 
Hence $1\cdot N_4=N_6\cdot{6\choose 4}$, and therefore $N_6=\frac{(q^4-1)(q^2-1)q^4}{15\cdot48}.$
\eop

%
%
\subsection{Proof of Theorem~\ref{Main thm 4}.}
Claim 1). Let $H=\{a^i,b^i,c^i\}_{i=1,2}$ be a centric hexagon of $\widehat{\X}$ with center $n_0$.

The lines $\langle a^1, b^1\rangle$ and $\langle a^2, b^2\rangle$ belong to the same plane hence they meet at a point $a^0\in H_0\cap Q$, where $H_0:=n_0^\perp.$
Similarly, the lines $\langle b^1, c^1\rangle$ and $\langle b^2, c^2\rangle$ meet at a point $b^0\in H_0\cap Q$ and the lines $\langle a^1, c^2\rangle$ and $\langle a^2, c^1\rangle$ meet at the point $c^0\in H_0\cap Q.$
Let $S:=\{a^i,b^i,c^i,\}_{i=0,1,2}$ and consider the subgeometry $Q(S)$ of $Q$
 induced on  $S$ (see Section~\ref{main results-sec 1}, paragraph before Theorem~\ref{Main thm 4}, for the definition), namely $Q(S):=(S, \mathcal{L}_S)$ where $\mathcal{L}_S=\{l\cap S\colon l {\rm \,\,is\,\,a
\,\,line\,\,of\,\,Q\,\,such\,\,that\,\, } |l\cap S|\geq 2\}.$

It is straightforward to see that the subgeometry $Q(S)$ is an $\F_2$-subgeometry isomorphic to a classical $(2,1)$-generalized quadrangle $Q^+(3,2)$. By construction, the points $a^0,b^0,c^0$ are orthogonal to $n_0$ (with respect to the orthogonality relation defined by $Q$) hence $H\cong Q^+(3,2)\setminus n_0^\perp.$

\medskip

Claim 2). Let $C=\{a^i,b^i,c^i,d^i\}_{i=1,2}$ be a cube of $\widehat{\X}$ with center $n_0$.
Let $\{\hat{l}_1, \hat{m}_1\}$, $\{\hat{l}_2, \hat{m}_2\}, \dots, \{\hat{l}_6, \hat{m}_6\}$ be the six pairs of opposite edges of $C$.
Then, for $j=1,\dots, 6$, the lines $\hat{l}_j$ and $\hat{m}_j$ have the same point at infinity $l_j^\infty=m_j^\infty=:p_j.$
These six points $p_1,p_2\dots, p_6$ are partitioned into three pairs and each of these pairs is collinear in $Q_0$ with a given point $p_0\in Q.$ Each of these three pairs of points (on a same line through $p_0$) corresponds to one of the three classes of parallel edges of the cube. 
The point $p_0$ is obtained as the sum of the $4$ vertices of any face of the cube.

Let $S:=\{a^i,b^i,c^i,d^i\}_{i=1,2}\cup \{p_i\}_{0\leq i\leq 6}$ and consider the subgeometry $Q(S)$ of $Q$
 induced on $S.$

It is straightforward to see that $Q(S)$ is an $\F_2$-subgeometry isomorphic to $Q(4,2)$.
It is also well known from that the complement in $Q(4,2)$ of a tangent hyperplane is a cube. Hence, $C\cong Q(S)\setminus p_0^\perp.$
Note that the projective subgeometry  $\langle S\rangle_{\F_2}$ defined by $S$ over the field $\F_2$ contains the nucleus $n_0$ of $Q_0$ but  $n_0$ is not the nucleus $n_S$ of the quadric $Q(S)$ (embedded in $\langle S\rangle_{\F_2}$).

\medskip

Claim 3) As remarked in Section~\ref{Introduction}, if $q=2$ then $\Gamma$ is a complete graph with $6$ vertices. In this case,  $Q\setminus Q_0$ is a dodecade.

Let $q\geq 4.$ Since by assumption $q$ is an odd power of $2$, we can take $\lambda=1$ in~(\ref{equation Q1}); thus $Q$ is represented by the following equation:
\begin{equation}\label{lambda=1}
x_1x_2+x_3x_4+x_5^2+x_5x_6+x_6^2=0.
\end{equation}

We shall count the number of subgeometries of $Q$ defined over the field $\F_2$ and represented by the equation~(\ref{lambda=1}) for a suitable choice of a basis.
Equivalently, we count the number of different bases $B_i$ of the vector space $V=V(6,q)$ with the property that with respect to each of these bases $B_i$, the quadratic form induced by $f$ on the $\F_2$-space spanned by $B_i$ is a quadratic form of elliptic type over $\F_2$, hence equivalent to $x_1x_2+x_3x_4+x_5^2+x_5x_6+x_6^2$.

Let us construct a basis  $B=(b_1,b_2,b_3,b_4,b_5,b_6)$ as required above. The first four vectors $b_1,b_2,b_3,b_4$ can be chosen so that $\langle b_1\rangle,$ $\langle b_2\rangle,$ $\langle b_3\rangle,$ $\langle b_4\rangle,$ is a quadrangle of $Q.$ That is, $f(b_i)=0$ for $1\leq i\leq 4$, $b_1\not\perp b_3$, $b_2\not\perp b_4$  and $b_1\perp b_2\perp b_3\perp b_4\perp b_1$.

The number of quadrangles of $Q$ is $(q^3+1)(q^2+1)(q+1)q^6.$

The condition  $\alpha(b_1,b_2)=1$ and $\alpha(b_3,b_4)=1$ (where $\alpha$ stands for  the linearization  of the form  $f$) still allows us to choose $b_1$ (or $b_2$) and $b_3$ (or $b_4$) up to arbitrary scalars. Hence,
we have $(q^3+1)(q^2+1)(q+1)q^6 (q-1)^2$  choices  for the quadruple $b_1,b_2,b_3,b_4.$  The vectors $b_5$ and $b_6$ are taken in the orthogonal space  $\langle b_1,b_2,b_3,b_4\rangle^{\perp}$ of $\langle b_1,b_2,b_3,b_4\rangle.$ (Recall that $\langle b_1,b_2,b_3,b_4\rangle^{\perp}$ is a projective line exterior to $Q.$)
So, let us regard $b_5$ and $b_6$ as vectors of the $2$-dimensional vector space $\langle b_1,b_2,b_3,b_4\rangle^{\perp}.$

Let $f'(x,y)=x^2+\mu xy+\nu y^2$ be the quadratic form induced by $f$ on $\langle b_1,b_2,b_3,b_4\rangle^{\perp},$ for a suitable choice of a basis $e_1$, $e_2$ of $\langle b_1,b_2,b_3,b_4\rangle^{\perp}.$   Hence $Tr(\nu/\mu^2)=1$ because
the line $\langle b_1,b_2,b_3,b_4\rangle^{\perp}$ is external to $Q.$  Put $b_5:=x_1 e_1 +x_2e_2$ and $b_6:=y_1e_1 +y_2e_2$  for scalars $x_1,x_2, y_1,y_2$ to be determined in a few lines.

The conditions to require on $b_5$ and $b_6$ are the following: $\alpha(b_5,b_6)=1$, $f(b_5)=1$ and $f(b_6)=1.$ By Lemma~\ref{lemma-algebraic condition}, we have $q+1$ different choices for $b_5.$ Indeed, $f(b_5)=1$ can be translated to $x_1^2+x_1x_2+x_2^2=1.$ 

Having chosen $b_5,$ let us count how many choices remain for $b_6.$

To fix ideas put $b_5=e_1.$ Then $\alpha(b_5,b_6)=1$ implies $y_2=1/\mu$ and  $f(b_6)=1$ implies $y_1^2+\mu y_1y_2+\nu y_2^2=1,$ namely $y_1^2+y_1+\nu/(\mu^2)+1=0.$
Since the equation $y_1^2+y_1+\nu/(\mu^2)+1=0$ has two distinct solutions, two different choices are left for $b_6.$ A similar argument works with any other choice of $b_5.$ Thus, in any case, only two choices are left for $b_6$ once $b_5$ has been chosen.
So, the number of $\F_2$-subgeometries $Q(S)$ of $Q^-(5,q)$ isomorphic to $ Q^-(5,2)$ is

\begin{equation}\label{lastequation}
{{ \begin{array}{l}
 \frac{\mbox{$\sharp$ choices of $(b_i)_{i=1}^{6}$ in $V(6,q)$}}{\mbox{$\sharp$ choices of  $(b_i)_{i=1}^{6}$ in  $V(6,2)$}}= \\
 \\
\frac{\mbox{$(q^3+1)(q^2+1)(q+1)q^6(q-1)^2(q+1)2$}}{\mbox{$(2^3+1)(2^2+1)(2+1)2^6(2-1)^2(2+1)2 $}}=\\
  \\
 \frac{\mbox{$(q^3+1)(q^2+1)(q+1)^2q^6(q-1)^2$ }}{\mbox{$9\cdot 5\cdot 3^2\cdot 2^6$}}.
 \end{array}}}
\end{equation}

By Proposition \ref{number of 6-cliques}, the total number $\bar{N}_6$ of centric dodecades of $Q$ (any possible center being allowed) is

\medskip

$ \begin{array}{l}
  \mbox{($\sharp$ of dodecades with a given  center)  $\cdot$  $|\PG(V)\setminus Q|=$}\\
   \\
\mbox{$ N_6\cdot |\PG(V)\setminus Q|=\frac{(q^4-1)(q^2-1)q^4}{15\cdot 48}\cdot \left[\frac{(q^6-1)}{q-1}-(q^3+1)(q+1))\right]=$}\\
 \\
\mbox{$\frac{(q^3+1)(q^4-1)(q^2-1)q^6}{48\cdot 15}$.}\\
 \end{array}$



\medskip

Hence $\bar{N}_6=\frac{(q^3+1)(q^4-1)(q^2-1)q^6}{48\cdot 15}.$
\medskip

Note that  if $Q(S)$ is an $\F_2$-subgeometry of $Q$ and $p\in \langle S\rangle_{\F_2}\setminus S$, then $Q(S)\setminus p^\perp$ is a dodecade with center $p.$ For each such subgeometry $Q(S)$ there are $2^2(2^3+1)=4\cdot 9$ possible choices for $p.$ In view of~(\ref{lastequation}) this accounts for a total of $\frac{(q^3+1)(q^2+1)(q+1)^2q^6(q-1)^2 }{9\cdot 5\cdot 3^2\cdot 2^6}\cdot 4\cdot 9=\frac{(q^3+1)(q^4-1)(q^2-1)q^6}{5\cdot9\cdot 2^4}=\bar{N}_6$ centric dodecades.

\medskip
Hence every dodecade with center $p$ corresponds to a subgeometry $Q(S)\setminus p^{\perp}$ over $\F_2$ as described above and viceversa.
\eop



\bigskip

\noindent
{\bf Acknowledgements.} The authors wish to thank Antonio Pasini for his very helpful remarks and comments on a first version of this paper.

\bigskip

\noindent
Authors' address\\

\noindent
Ilaria Cardinali,  \\
Department of Information Engineering and Mathematics,\\
University of Siena,\\
Via Roma 56, 53100 Siena, Italy\\
ilaria.cardinali@unisi.it,

\bigskip

\noindent
Narasimha N. Sastry,  \\
Theoretical Statistics and Mathematics Unit \\
Indian Statistical Institute,\\
8th Mile Mysore Road,\\
560059 Bangalore, India\\
nsastry@isibang.ac.in.

\end{document}